\documentclass{amsproc}

\usepackage[dvips]{graphicx}
\usepackage{rotating}
\usepackage{paralist}
\usepackage{psfrag}
\usepackage{lscape}
\usepackage{enumerate}
\usepackage{url}
\usepackage{subfigure}

\newcommand{\R}{\ensuremath{\mathbb{R}}}
\newcommand{\polar}[1]{\ensuremath{{#1}^{\Delta}}}
\newcommand{\fcone}{\ensuremath{\mathcal D}}
\newcommand{\fhull}{\ensuremath{\mathcal C}}
\newcommand{\skp}[2]{\langle #1,#2\rangle}
\newcommand{\pstack}[3]{\mathbf{PS}_{#2}^{#1}(#3)}
\newcommand{\stdpstack}{\pstack{S}{\mathcal F, \mathcal N}{P}}
\newcommand{\adj}[1]{\mathrm{adj}(#1)}
\newcommand{\hadj}[1]{\mathcal H_{#1}}
\newcommand{\padj}[1]{\mathcal H^+_{#1}}
\newcommand{\nv}[3]{\mathcal R^#1_{#2}(#3)}
\DeclareMathOperator{\facets}{Fac}
\newcommand{\hfacets}{{\mathcal H}}
\newcommand{\pfacets}{{\mathcal H^+}}
\newcommand{\first}{\mathcal I^1}
\newcommand{\second}{\mathcal I^2}
\newcommand{\eps}{\varepsilon}
\newcommand{\tsts}{$2$-simplicial, $2$-simple}
\newcommand{\tstsfp}{\tsts\ $4$-polytope}
\DeclareMathOperator{\fat}{F}
\DeclareMathOperator{\compl}{C}
\DeclareMathOperator{\conv}{conv}
\DeclareMathOperator{\flag}{flag}
\DeclareMathOperator{\fdeg}{fdeg}
\newenvironment{verticestable}{\begin{displaymath}\begin{array}{@{\vspace{2mm}\big(}r@{,\;}r@{,\;}r@{,\;}r@{\big)}}}{\end{array}\end{displaymath}}
\newenvironment{viftable}{\begin{displaymath}\begin{array}{@{\vspace{2mm}}l}}{\end{array}\end{displaymath}}

\theoremstyle{plain}
\newtheorem{theorem}{Theorem}[section]
\newtheorem{corollary}[theorem]{Corollary}
\newtheorem{lemma}[theorem]{Lemma}
\newtheorem{proposition}[theorem]{Proposition}

\newtheorem*{theorem*}{Theorem \ref{thm:many-polys}}
\newtheorem*{corollary*}{Corollary \ref{cor:ray}}

\theoremstyle{definition}
\newtheorem{definition}[theorem]{Definition}

\theoremstyle{remark}
\newtheorem{remark}[theorem]{Remark}

\begin{document}

\title[Constructions for $4$-Polytopes and the Cone of Flag Vectors]
      {Constructions for $4$-Polytopes \\ and the Cone of Flag Vectors}

\author{Andreas Paffenholz}
\address{FU Berlin\\ 
         Institut f\"ur Mathematik II\\ 
         Arnimallee 3\\ 
         14195 Berlin\\ 
         Germany}

\email{paffenho@math.tu-berlin.de}

\thanks{The   first     author  was      supported by    the  Deutsche
  Forschungsgemeinschaft    within  the    European   graduate program
  `Combinatorics, Geometry, and Computation' (No. GRK 588/2)}

\author{Axel Werner}
\address{TU Berlin\\ 
         Institut f\"ur Mathematik\\ 
         MA 6-2\\ 
         Stra{\ss}e des 17.~Juni 136\\ 
         10623 Berlin\\ 
         Germany}

\email{awerner@math.tu-berlin.de}

\thanks{The second author was supported by the DFG Leibniz grant of G.~M.~Ziegler.}


\date{November 30, 2005}


\begin{abstract}
  We describe a  construction for $d$-polytopes generalising the  well
  known stacking operation.  The   construction is applied to  produce
  $2$-simplicial  and $2$-simple   $4$-polytopes with $g_2=0$   on any
  number of $n \geq 13$ vertices. In particular, this implies that the
  ray $\ell_1$,  described by Bayer  (1987), is fully contained in the
  convex hull  of all   flag  vectors of   $4$-polytopes.   Especially
  interesting examples on $9$, $10$ and $11$ vertices are presented.
\end{abstract}

\maketitle

%%%%%%%%%%%%%%%%%%%%%%%%%%%%%%%%%%%%%%%%%%%%%%%%%%%%%%%%%%%%%%%%%%%%%%%%%%%%%%%

\section{Introduction} \label{sec:intro}

It  is a wide  open problem  in  discrete  geometry to understand  the
combinatorial properties  of polytopes, which  can be described as the
convex hull of finitely  many points  in  some $\R^n$.  Even the  more
special task  to characterise the $f$- and  flag vectors  of polytopes
(and more  generally of  spheres)  of arbitrary dimension  seems to be
very hard.   For  $3$-dimensional  polytopes, Steinitz \cite{steinitz}
gave  a   complete  characterisation.   For arbitrary  dimension,  the
problem is  still  open,  although some   conditions  are known.    In
particular,  all linear relations between  the entries of flag vectors
are    described   by  the  Generalized  Dehn-Sommerville    equations
\cite{MR86f:52010b}  and the  admissible   $f$-vectors  of  simplicial
polytopes are classified in terms of their $g$-vectors
by the $g$-theorem of Billera, Lee and Stanley and McMullen
(see for instance \cite[Thm.~8.35]{MR96a:52011}). Additionally, for dimension $4$ a
linear  approximation of the set of  flag vectors  and $f$-vectors was
given by Bayer \cite{MR88b:52009}.  A different view on the $f$-vector
characterisation  was provided   by  Ziegler  \cite{MR1957565}.  Since
then,   some   progress   has    been    made   (cf.~\cite{MR2096750},
\cite{MR2119033}), but it is still not known what the linear cone of
$f$- resp.~flag vectors looks like.

In this paper we introduce a new polytope construction method,
examine some of its  combinatorial properties and construct $4$-polytopes with
arbitrarily high numbers of vertices   that are extremal for the  flag
vector cone.  In particular, we  prove the following two theorems (see
below for definitions).
\begin{theorem*}
  Elementary \tstsfp{s} with $k$ vertices exist for $k=5,9,10,11$ and
  $k \geq 13$.
\end{theorem*}
This implies (using the notation of Bayer \cite{MR88b:52009}):
\begin{corollary*}
  The ray $\ell_1$ is contained in the convex hull
  of all flag vectors of $4$-polytopes.
\end{corollary*}
Additionally,   we briefly analyse  the   consequences of the  various
recent polytope constructions in \cite{MR2096750} and \cite{MR2119033}
for the flag vector cone.

\subsection*{Acknowledgements}

The authors would like to thank G\"unter M.\ Ziegler, who 
supported this work at several occasions and in particular
pointed out the significance of our series of polytopes,
Eran Nevo and Raman Sanyal for helpful discussions, 
and the referee for several suggestions improving the exposition.

%%%%%%%%%%%%%%%%%%%%%%%%%%%%%%%%%%%%%%%%%%%%%%%%%%

\subsection{General preliminaries}

We first give the basic definitions as well as an overview over
important related results.

Let $P$  be   a $d$-polytope and $[d]  =   \{ 0,\ldots,d-1 \}$.    The
\emph{$f$-vector} of $P$ is the $d$-dimensional vector
\begin{displaymath}
  f(P) \; = \; (f_0,f_1,\ldots,f_{d-1}) ,
\end{displaymath}
where $f_i$   for  $0  \leq  i    \leq  d-1$ denotes the   number   of
$i$-dimensional faces of $P$.   The \emph{flag vector}  of $P$ is  the
$2^d$-dimensional vector
\begin{displaymath}
  \flag(P) \; = \; (f_S)_{S \subseteq [d]} ,
\end{displaymath}
where $f_S$ for  a subset $S =  \{ i_1,\ldots,i_k \}$ of $[d]$ denotes
the number of face  chains $\emptyset \subset F_{i_1} \subset \ldots
\subset  F_{i_k} \subset P$ such  that  $\dim F_{i_j} =  i_j$ for $1
\leq j \leq k$.  We usually write  $f_{i_1 i_2 \ldots i_k}$ instead of
$f_{\{i_1,i_2,\ldots,i_k\}}$. Faces of  codimension  $1$, $2$ and  $3$
are   called   \emph{facets},    \emph{ridges}    and \emph{subridges}
respectively.

The   Generalized  Dehn-Sommerville  equations   by Bayer and  Billera
\cite{MR86f:52010b} imply that the flag vectors of $d$-polytopes lie in
an   $(F_d-1)$-dimensional affine   subspace    of  $\R^{2^d}$,  where
$(F_k)_{k \geq  0}  = (1,1,2,3,\ldots)$  is   the series of  Fibonacci
numbers.

$P$ is \emph{$k$-simplicial} ($0 \leq k \leq d-1$) if all its $k$-dimensional
faces  are  $k$-simplices,   i.e.~contain   exactly   $k+1$  vertices;
equivalently,  if    in  the face   lattice   of  $P$  all   intervals
$[\emptyset,F]$ for $\dim F = k$ are boolean.

$P$ is \emph{$h$-simple} ($0 \leq h \leq d-1$) if the dual polytope
$\polar{P}$ is  $h$-simplicial; equivalently, if all intervals $[F,P]$
with  $\dim  F  = d-h-1$ are    boolean or equivalently,  if  all  its
$(d-h-1)$-dimensional faces are contained in exactly $h+1$ facets.

The \emph{$g$-vector} of a polytope can be defined in general
using generating  functions (see  \cite{MR98a:05001}); its entries can
be written as linear  combinations of the  entries of the flag vector.
We will, however, focus on one special entry.
\begin{definition}
  For a $d$-polytope $P$ define
  \begin{align}\label{g2}
    g_2(P) \; = \; f_{02} - 3 f_2 + f_1 - d f_0 + {d+1 \choose 2} .
  \end{align}
\end{definition}
It has been shown by Kalai \cite{MR877009} via rigidity theory
that $g_2(P) \geq 0$  for every $d$-polytope $P$ with $d \geq 4$.
Polytopes $P$ with $g_2(P)=0$ are called \emph{elementary}.
It is an interesting open problem to characterise all elementary polytopes.
See also \cite{Kalai94} for a survey on this topic.

%%%%%%%%%%%%%%%%%%%%%%%%%%%%%%%%%%%%%%%%%%%%%%%%%%

\subsection{Flag vectors of $4$-polytopes}

In the last  two sections  we will focus  on  $4$-polytopes, where the
situation is  understood slightly better than  in general.  It is easy
to prove that for \tstsfp{s} the $f$-vector is  symmetric and the flag
vector is completely determined by the $f$-vector.

\begin{figure}[t]
  \centering
  \begin{picture}(0,0)%
    \includegraphics{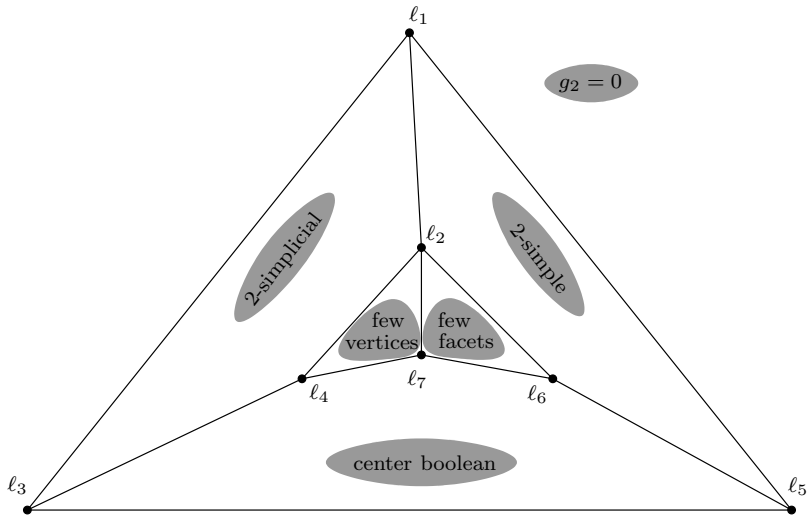}%
  \end{picture}%
  \setlength{\unitlength}{3947sp}%
  \begingroup\makeatletter\ifx\SetFigFont\undefined%
  \gdef\SetFigFont#1#2#3#4#5{%
    \reset@font\fontsize{#1}{#2pt}%
    \fontfamily{#3}\fontseries{#4}\fontshape{#5}%
    \selectfont}%
  \fi\endgroup%
  \begin{picture}(4957,3202)(476,-2792)
    \put(3941,-106){\makebox(0,0)[lb]{\smash{{\SetFigFont{8}{9.6}{\ttdefault}{\mddefault}{\updefault}{$g_2=0$}%
          }}}}
    \put(3001,314){\makebox(0,0)[lb]{\smash{{\SetFigFont{8}{9.6}{\ttdefault}{\mddefault}{\updefault}{$\ell_1$}%
          }}}}
    \put(3107,-1062){\makebox(0,0)[lb]{\smash{{\SetFigFont{8}{9.6}{\ttdefault}{\mddefault}{\updefault}{$\ell_2$}%
          }}}}
    \put(3609,-1001){\makebox(0,0)[lb]{\smash{{\SetFigFont{8}{9.6}{\familydefault}{\mddefault}{\updefault}{\begin{turn}{308}$2$-simple\end{turn}}%
          }}}}
    \put(1951,-1486){\makebox(0,0)[lb]{\smash{{\SetFigFont{8}{9.6}{\familydefault}{\mddefault}{\updefault}{\begin{turn}{52}$2$-simplicial\end{turn}}%
          }}}}
    \put(2765,-1607){\makebox(0,0)[lb]{\smash{{\SetFigFont{8}{9.6}{\familydefault}{\mddefault}{\updefault}{few}%
          }}}}
    \put(3180,-1607){\makebox(0,0)[lb]{\smash{{\SetFigFont{8}{9.6}{\familydefault}{\mddefault}{\updefault}{few}%
          }}}}
    \put(3185,-1740){\makebox(0,0)[lb]{\smash{{\SetFigFont{8}{9.6}{\familydefault}{\mddefault}{\updefault}{facets}%
          }}}}
    \put(2599,-1746){\makebox(0,0)[lb]{\smash{{\SetFigFont{8}{9.6}{\familydefault}{\mddefault}{\updefault}{vertices}%
          }}}}
    \put(2651,-2505){\makebox(0,0)[lb]{\smash{{\SetFigFont{8}{9.6}{\familydefault}{\mddefault}{\updefault}{center boolean}%
          }}}}
    \put(476,-2642){\makebox(0,0)[lb]{\smash{{\SetFigFont{8}{9.6}{\ttdefault}{\mddefault}{\updefault}{$\ell_3$}%
          }}}}
    \put(2988,-1974){\makebox(0,0)[lb]{\smash{{\SetFigFont{8}{9.6}{\ttdefault}{\mddefault}{\updefault}{$\ell_7$}%
          }}}}
    \put(5381,-2648){\makebox(0,0)[lb]{\smash{{\SetFigFont{8}{9.6}{\ttdefault}{\mddefault}{\updefault}{$\ell_5$}%
          }}}}
    \put(3743,-2068){\makebox(0,0)[lb]{\smash{{\SetFigFont{8}{9.6}{\ttdefault}{\mddefault}{\updefault}{$\ell_6$}%
          }}}}
    \put(2374,-2061){\makebox(0,0)[lb]{\smash{{\SetFigFont{8}{9.6}{\ttdefault}{\mddefault}{\updefault}{$\ell_4$}%
          }}}}
  \end{picture}%
  \caption{Hyperplane section through the cone \fcone, according to \cite{MR88b:52009}}
  \label{fig:flagcone-dim4-bayer}
\end{figure}

Bayer \cite{MR88b:52009} described all the known linear inequalities for
flag vectors of $4$-polytopes.  Let $\fcone \subseteq \R^4$ be the polyhedron
defined by these inequalities and $\fhull \subseteq \fcone$ the convex hull
of all flag vectors of $4$-polytopes.
Then \fcone\ is a $4$-dimensional cone with the flag vector of the $4$-simplex as its apex.
Figure~\ref{fig:flagcone-dim4-bayer} illustrates a hyperplane section through this cone.
It can be viewed as a $3$-polytope with its  vertices, edges and faces
representing  special properties of $4$-polytopes,  in  the sense that
polytopes whose flag  vectors lie on  the respective faces  have these
properties.
The main question is how close the cone \fcone\ approximates
\fhull. The task is therefore to  find examples of polytopes with flag
vectors in extremal regions of \fcone.

In Table~\ref{tab:rays-of-flagcone}
we give a summary of  what is known for the  rays of \fcone; note that
$\ell_4$ and $\ell_6$ are not contained  in \fhull\ itself, but in its
closure (cf.~\cite[Sec.~2]{MR88b:52009}).
\begin{table}[b]
  \caption{Known polytopes on or close to the rays of \fcone}
  \label{tab:rays-of-flagcone}
  \begin{tabular}{cll}
    ray & property & examples \\ \hline
    $\ell_1$ & \tstsfp{s} with $g_2=0$ & see Theorem~\ref{thm:many-polys} \\
    $\ell_2$ & `fat' \tstsfp{s} & unknown \\
    $\ell_7$ & `fat' center boolean $4$-polytopes & unknown \\
    $\ell_4$ & simplicial $4$-polytopes with few vertices & cyclic polytopes \\
    $\ell_3$ & simplicial $4$-polytopes with $g_2=0$ & stacked polytopes \\
    $\ell_5$ & simple $4$-polytopes with $g_2=0$ & truncated polytopes \\
    $\ell_6$ & simple $4$-polytopes with few facets & dual cyclic polytopes
  \end{tabular}
\end{table}

Until now only two flag vectors of polytopes on $\ell_1$ were known,
that of the $4$-simplex $\Delta_4$ and that of the $4$-dimensional hypersimplex
$\Delta_4(2)$  (for a   definition  see \cite[Ch.~0]{MR96a:52011}   for
instance).  Hence  it was  clear  that $\ell_1$  contained an edge  of
\fhull,  but not,  whether $\ell_1$ is a ray of \fhull, nor if $\ell_1$ 
contained any further flag vectors at all.
We  establish in   Section~\ref{sec:ex} that   $\ell_1$ is
indeed an extremal ray of $\fhull$.

A few more regions of \fcone\ deserve a closer study. The set of flag
vectors of general \tstsfp{s} is a subset of the $2$-dimensional cone spanned
by $\ell_1$ and $\ell_2$. It is indeed  a $2$-dimensional set -- there
are \tstsfp{s} with the same number of vertices, but different numbers
of edges (cf.~\cite{MR2096750}).  Analogously, the set of flag vectors
of elementary $4$-polytopes  is contained in the $3$-dimensional  cone
spanned by $\ell_1$, $\ell_3$ and $\ell_5$. There are also a number of
extreme examples: stacked polytopes,  iterated pyramids  over $n$-gons
and  multiplexes   (cf.~\cite{MR1871688}).  Also,  the   constructions
described in the next sections produce in general elementary polytopes
when    applied   to       such,     as   can       be   seen     from
Corollary~\ref{cor:g2}. Note that $2$-simplicity,
$2$-simpliciality and  $g_2=0$ are three independent  properties, that
is, there are $4$-polytopes with any combination of these properties.

Bayer conjectured that the hyperplane determined by $\ell_2$, $\ell_4$
and $\ell_6$ yields a valid inequality for $4$-polytopes.  This is not
true, since there exist polytopes with flag vectors close to the interior of the edges
$[\ell_4,\ell_7]$, and dually $[\ell_6,\ell_7]$, as shown by Joswig \&
Ziegler~\cite{MR1758054} and   Ziegler~\cite{MR2119033}.    However, a
hyperplane cutting off  $\ell_7$  may  still be possible.
Extremal polytopes known in this respect are projected products of polygons
by   Ziegler~\cite{MR2119033}    and   polytopes   obtained   by   the
$E$-construction  (Paffenholz~\& Ziegler~\cite{MR2096750}).  The  most
restrictive  linear inequality that  is compatible with these examples
would then be
\begin{displaymath}
  f_{03} - 140 \; \geq \; 4 (f_1+f_2) - 20 (f_0+f_3) .
\end{displaymath}
In terms of the  \emph{fatness} $\fat$ and \emph{complexity}  $\compl$
of $4$-polytopes,  as   introduced by Ziegler  \cite{MR1957565},  this
inequality reads $4\fat -\compl\le 20$.  In this respect, providing an
upper  bound for fatness, that  is,  bounding the  number of edges and
ridges by the number of vertices and facets, would be helpful.

%%%%%%%%%%%%%%%%%%%%%%%%%%%%%%%%%%%%%%%%%%%%%%%%%%%%%%%%%%%%%%%%%%%%%%%%%%%%%%%

\section{General Construction} \label{sec:constr-gen}

This section consists of two parts.  In the first  part we provide the
basic tool for the construction of a family of $2$-simple and
$2$-simplicial $4$-polytopes with vanishing $g_2$. The operation
has many more applications than the ones we will discuss in more detail in
the  second part. Hence,  we  give all   definitions and theorems  for
arbitrary  dimension  $d\ge  3$. We  indicate   some of the additional
applications we have encountered so far.

In the second part we examine settings in which the operation can
be  applied,  determine  the  facet  types   that  can occur  in   the
construction and prove  that some interesting  properties of polytopes
are preserved.  This part will  be much more specifically tailored for
what we  need for the polytope families  defined later, as some of the
$f$-vector computations  tend to get  complicated  in a  more  general
setting.

In Section~\ref{sec:constr-2s2s}  we combine several  instances of our
basic tool to obtain two special  constructions $\first$ and $\second$
producing   $4$-polytopes  $\first(P;S)$ and   $\second(P;S)$ out of a
polytope $P$ and a facet  $S$  of $P$.   Using  the properties of  the
construction discussed  in this section  we prove that $2$-simplicity,
$2$-simpliciality  and  the  value   of  $g_2$  is   preserved by  the
construction.

%%%%%%%%%%%%%%%%%%%%%%%%%%%%%%%%%%%%%%%%%%%%%%%%%%

\subsection{Pseudo-Stacking}

The operation we introduce here is a
generalisation of the well known stacking operation.  In most cases it
adds one new vertex to the polytope. We need some notation for this.
\begin{definition}
  Let $P$ be a $d$-polytope. A \emph{simplex facet} of $P$ is a facet
  of $P$ that is combinatorially equivalent to a $(d-1)$-simplex.
\end{definition}
Let $H:=\{x\in\R^d\mid   \skp{x}{v}=\ell\}$ for some $v\in   \R^d$ and
$\ell\in\R$   be   an   affine hyperplane.    By  $H^+:=\{x\in\R^d\mid
\skp{x}{v}>\ell\}$  we denote  the positive half  space  defined  by $H$
and similarly the negative half space by $H^-$.

Let   $F$ be a facet  of  a $d$-polytope $P$.  We  denote by $H_F$ the
unique affine hyperplane  that contains $F$,  oriented  in such  a way
that  $P$ is contained in  $H^+_F\cup H^{\phantom{+}}_F$.  By $\facets(P)$ we denote
the set of all facets of $P$,  by $\hfacets(P) := \{H_F \mid F \in \facets(P)\}$
the set of all hyperplanes coming from the  facets, and  by $\pfacets(P)$ the set
$\{H_F^+\mid F \in \facets(P)\}$ of half-spaces.

Further let  $\adj{F}$ be the set  of  facets of $P$  adjacent to $F$.
For  a  subset  $\mathcal  F$ of  $\adj{F}$   we  denote  the set   of
hyperplanes defined  by the facets  in $\mathcal F$ by $\hadj{\mathcal
  F}$,  and the   set   of  positive  half-spaces determined   by  the
hyperplanes in $\hadj{\mathcal F}$ by $\padj{\mathcal F}$.

Finally, if $v$ is some point in $\R^d$  then $v$ lies \emph{beyond} a
facet $F$ if $v  \in H^-_F$ and $v$ lies  \emph{beneath} $F$ if $v \in
H^+_F$.

Let $S$ be a facet of  a $d$-polytope $P$  and $\mathcal F,\mathcal N$
be disjoint subsets of $\adj{S}$. We define the region
$\nv{S}{\mathcal F, \mathcal N}{P}$ in $\R^d$ by
\begin{align*}
  \nv{S}{\mathcal F, \mathcal N}{P} & :=
  \left( \bigcap_{H \in \mathcal A} H^+ \right) \cap
  \left( \bigcap_{H \in \hadj{\mathcal F}} H \right) \cap
  \left( \bigcap_{H \in \hadj{\mathcal N} \cup \{H_S\}} H^- \right)
\end{align*}
for $\mathcal A := \hfacets(P) \setminus ( \hadj{\mathcal F} \cup \hadj{\mathcal N} \cup \{H_S\} )$.

\begin{definition}[Pseudo-Stacking]
  Let $P$  be a $d$-polytope and  $S$ a simplex facet  of $P$.  Choose
  two disjoint sets $\mathcal F, \mathcal N \subseteq \adj{S}$.
  Assume that $\nv{S}{\mathcal F,\mathcal N}{P} \ne \emptyset$.

  The \emph{pseudo-stacking $\pstack{S}{\mathcal F,\mathcal N}{P}$  of
    $P$ above $S$ with  respect to $\mathcal F$  and $\mathcal  N$} is
  the   convex    hull   $\conv(P\cup  v)$  of    $P$    with a  point
  $v\in\nv{S}{\mathcal F,  \mathcal N}{P}$.  If $\mathcal N=\emptyset$
  then we omit it in    the notation and write    $\pstack{S}{\mathcal
    F}{P}$ for the pseudo-stacking of $P$ above $S$.
\end{definition}
Compare this definition to Gr\"un\-baum \cite[Section~5.2]{MR1976856}.
Also, a similar concept was studied by Altshuler and Shemer \cite{MR738161}
for the purpose of enumeration of $4$-polytopes with few vertices.

\begin{figure}[b]
  \centering
  \psfrag{v0}[t][t]{$v_0$}
  \psfrag{v1}[l][l]{$v_1$}
  \psfrag{v2}[bl][bl]{$v_2$}
  \psfrag{v3}[r][r]{$v_3$}
  \psfrag{v4}[r][r]{$v_4$}
  \psfrag{v5}[bl][bl]{$v_5$}
  \psfrag{v6}[b][b]{$v_6$}
  \psfrag{v7}[l][l]{$v_7$}
  \psfrag{S}[r][r]{$S$}
  \psfrag{F1}[r][r]{$F_1$}
  \psfrag{F2}[l][l]{$F_2$}
  \psfrag{F3}[l][l]{$F_3$}
  \begin{minipage}[t]{.49\textwidth}
    \includegraphics[width=.99\textwidth]{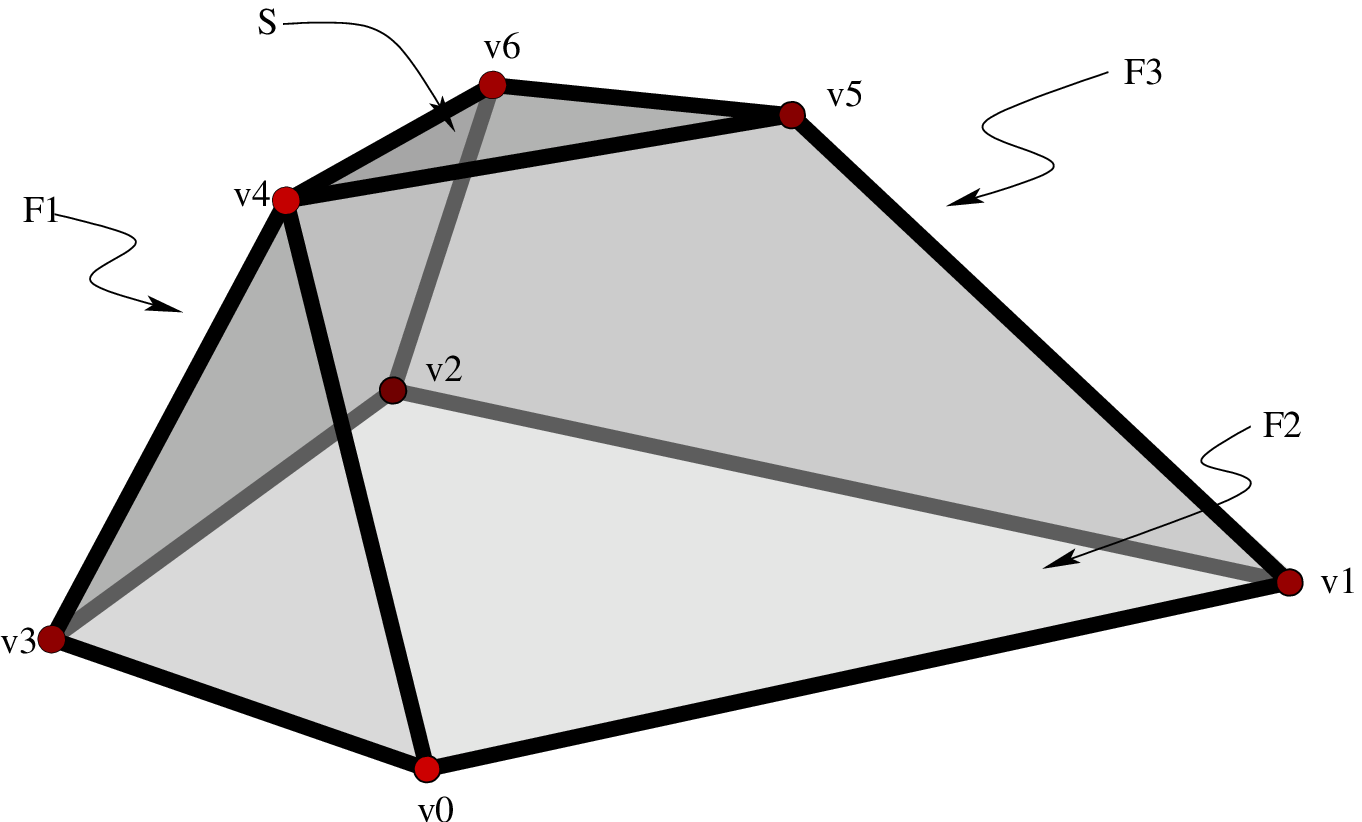}
    \caption{A polytope $P$ with simplex facet $S$ in bounded position
      and $\adj{S}=\{F_1,F_2,F_3\}$.  }
    \label{fig:ps01}
  \end{minipage}
  \hspace{.01\textwidth}
  \begin{minipage}[t]{.48\textwidth}
  \psfrag{v6}[br][br]{$v_6$}
    \includegraphics[width=.99\textwidth]{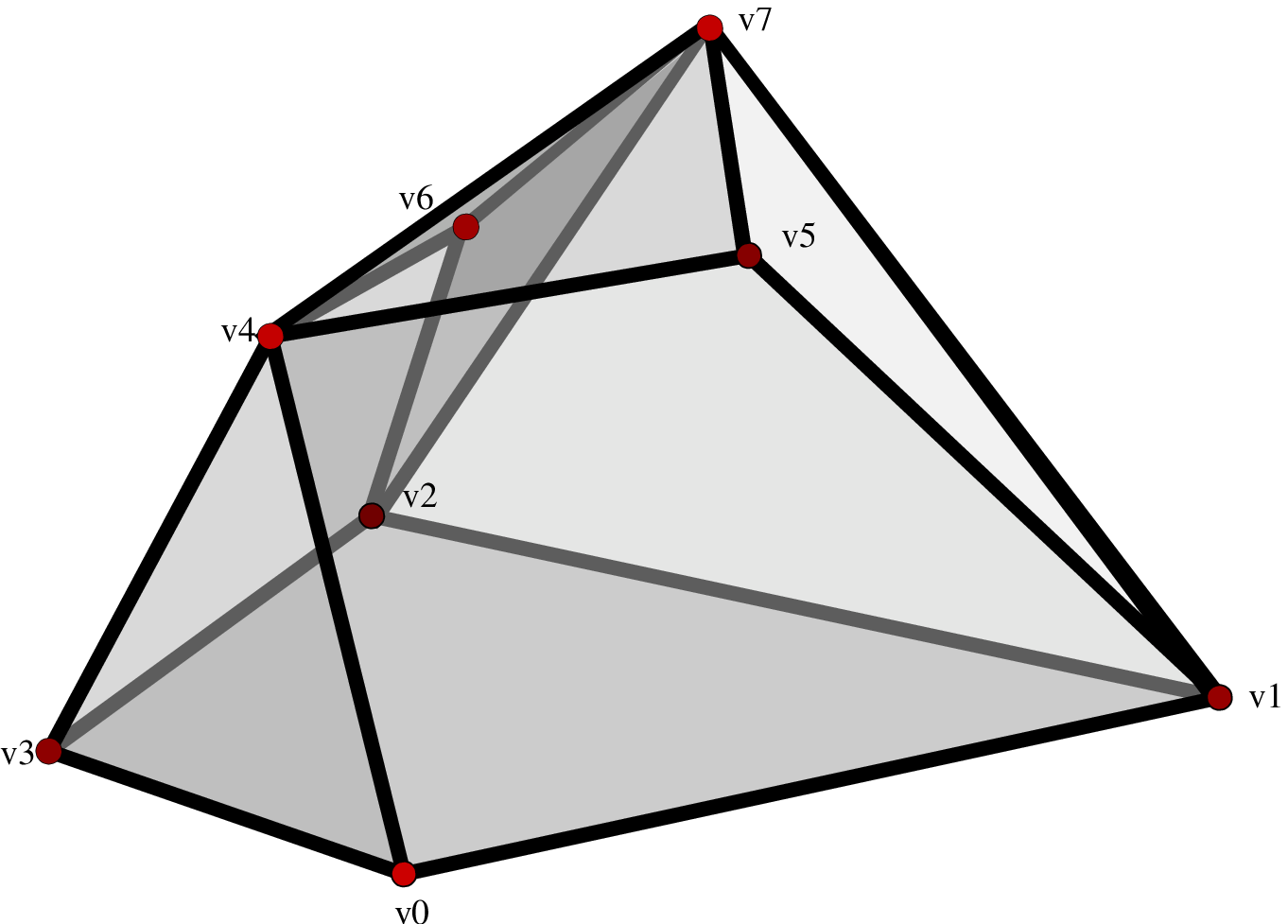}
    \caption{The polytope $\pstack{S}{\mathcal F, \mathcal N}{P}$ with
      $\mathcal F:=\emptyset$ and $\mathcal N:=\{F_3\}$.}
    \label{fig:ps02}
  \end{minipage}
\end{figure}

\begin{figure}[b]
  \centering
  \psfrag{v0}[t][t]{$v_0$}
  \psfrag{v1}[l][l]{$v_1$}
  \psfrag{v2}[bl][bl]{$v_2$}
  \psfrag{v3}[r][r]{$v_3$}
  \psfrag{v4}[r][r]{$v_4$}
  \psfrag{v5}[bl][bl]{$v_5$}
  \psfrag{v6}[b][b]{$v_6$}
  \psfrag{v7}[l][l]{$v_7$}
\begin{minipage}[t]{.48\textwidth}
  \psfrag{v6}[br][br]{$v_6$}
    \includegraphics[width=.99\textwidth]{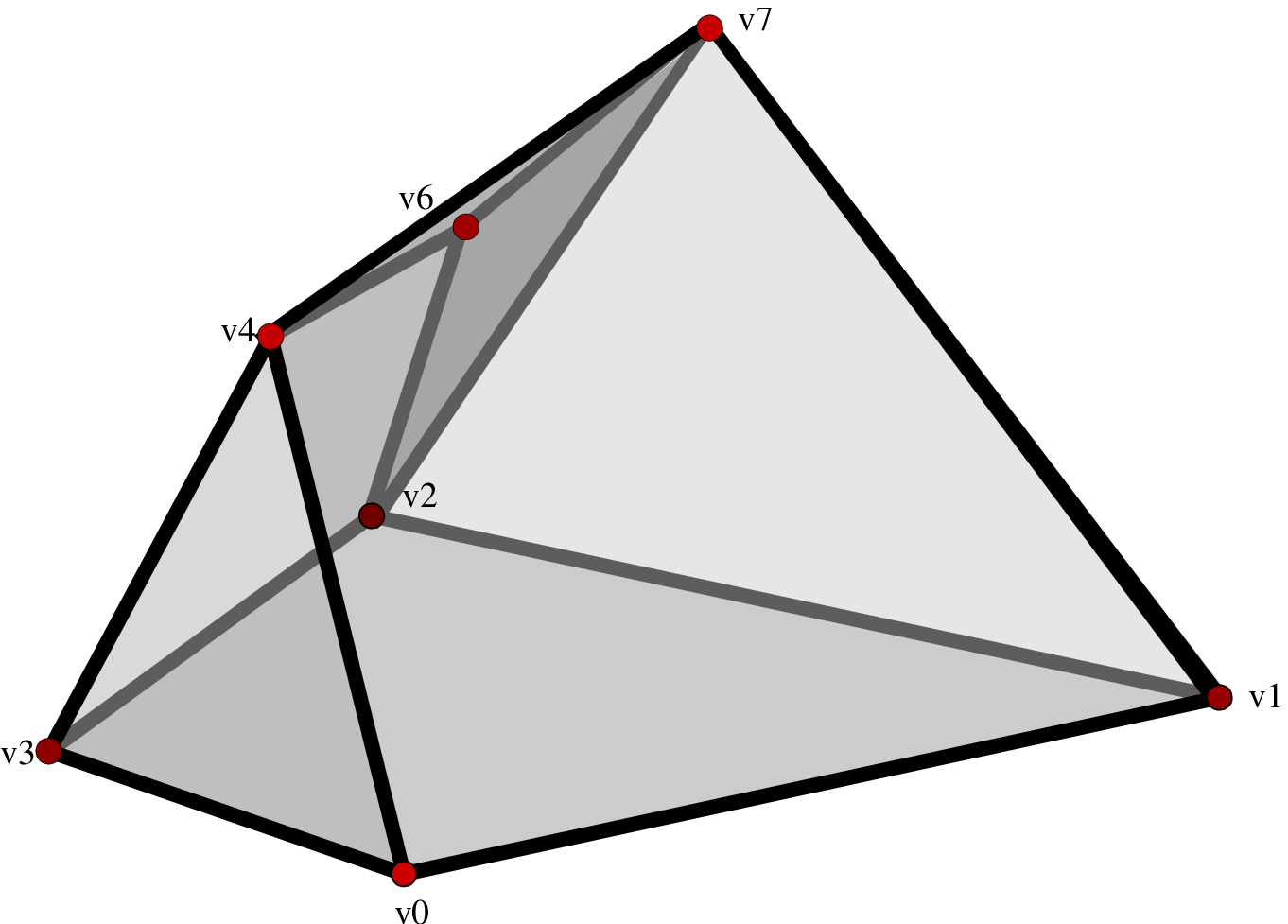}
    \caption{The polytope $\pstack{S}{\mathcal F, \mathcal N}{P}$ with
      $\mathcal F:=\{F_2\}$, $\mathcal N:=\{F_3\}$.}
    \label{fig:ps03}
  \end{minipage}
  \hspace{.02\textwidth}
  \begin{minipage}[t]{.48\textwidth}
    \psfrag{v6}[br][br]{$v_6$}
    \includegraphics[width=.99\textwidth]{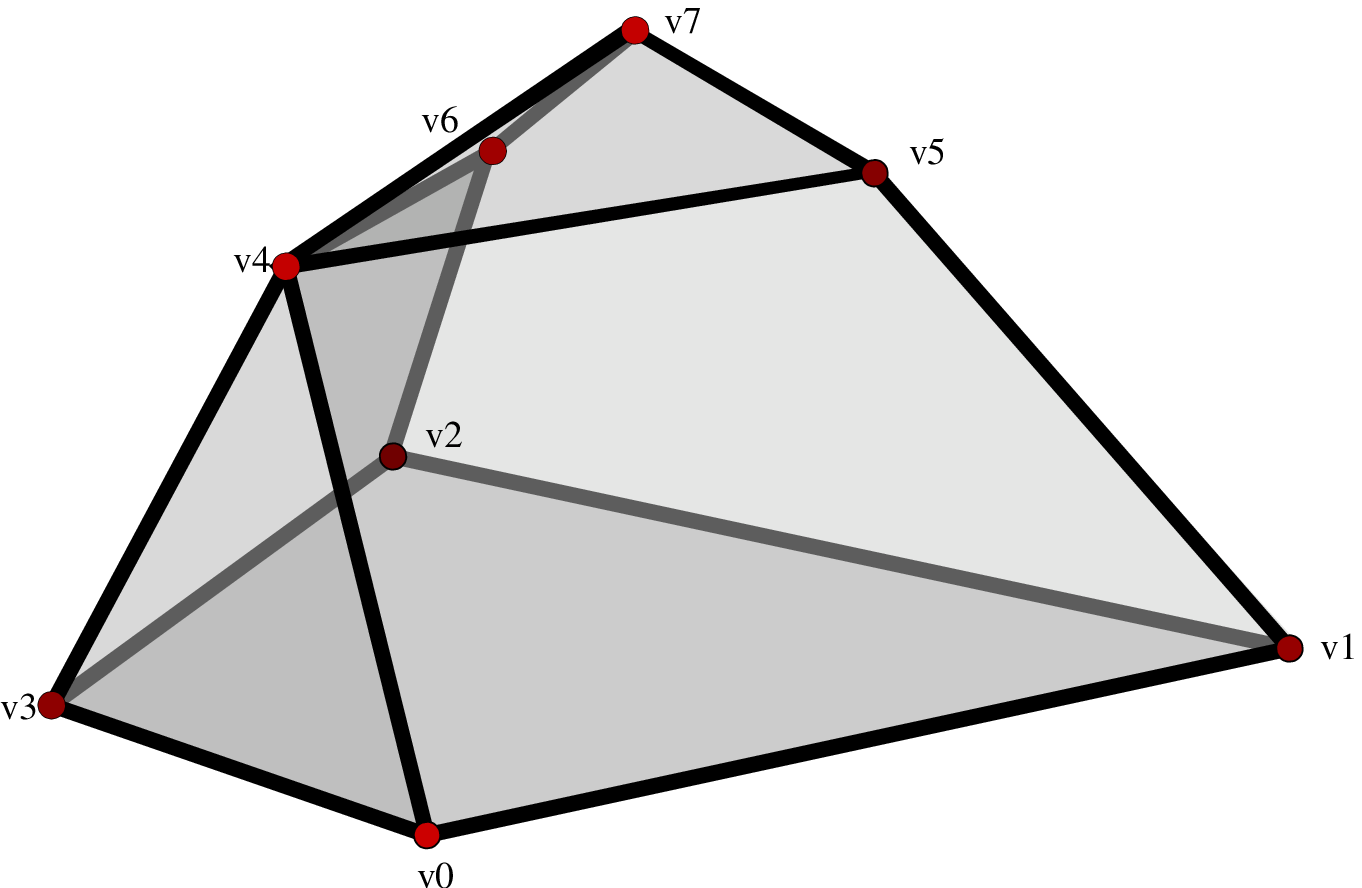}
    \caption{The polytope $\pstack{S}{\mathcal F, \mathcal N}{P}$ with
      $\mathcal F:=\{F_3\}$ and $\mathcal N:=\emptyset$.}
    \label{fig:ps04}
  \end{minipage}
\end{figure}

In plain words, a point $v\in\nv{S}{\mathcal F, \mathcal N}{P}$ lies
\textit{beyond} $S$ and all facet hyperplanes coming from facets in $\mathcal N$,
it lies \textit{in} all facet hyperplanes from facets in $\mathcal F$,
and \textit{beneath} all other facet hyperplanes. A priori, the set $\nv{S}{\mathcal F, \mathcal N}{P}$ need
not contain a point at all, but  we will show  some conditions that guarantee
that this set is non-empty.  See Figures~\ref{fig:ps01}--\ref{fig:ps04}
for some illustrations  of   this operation with various   choices  of
$\mathcal F$ and $\mathcal N$.
We have the following simple fact about the pseudo-stacking operation:
\begin{proposition}
  The combinatorial   properties   of $\pstack{S}{\mathcal  F,\mathcal
    N}{P}$    do not  depend   on the   actual   choice  of  the point
  $v\in\nv{S}{\mathcal F, \mathcal N}{P}$.\qed
\end{proposition}

\begin{remark}
  The   usual stacking  operation is   the  special case  $\mathcal F=
  \mathcal N= \emptyset$.
\end{remark}

%%%%%%%%%%%%%%%%%%%%%%%%%%%%%%%%%%%%%%%%%%%%%%%%%%

\subsection{Properties of the Construction}

Now we want to examine some cases in which the set
$\nv{S}{\mathcal F, \mathcal N}{P}$ is non-empty and thus
pseudo-stacking can be applied. To simplify the statements of the
propositions we introduce some more notation.
\begin{definition}
  Let $P$ be  a $d$-polytope   and $F$  a facet   of  $P$.  A   subset
  $\mathcal F$ of $\adj{F}$ is called  \emph{nonsimple} if there is no
  pair   of  adjacent facets $G,G'\in\mathcal  F$   that have a common
  $(d-3)$-face with $F$.
\end{definition}
Nonsimplicity of the set $\mathcal F$  implies that there is no simple
$(d-3)$-face (i.e.\ a $(d-3)$-face  contained in precisely $3$ facets)
which is a  subface of $S$ and  the other two  facets are in $\mathcal
F$.
\begin{definition}
  Let $P$ be a  $d$-polytope. A facet $F$  of $P$ is in  \emph{bounded
    position}   if  the hyperplanes  in any   subset of $\hadj{\adj{F}}$ of
  cardinality $d$ intersect in  a point in $H_F^-$.  Equivalently, $F$
  is  in bounded   position  if the intersection    of the half-spaces
  defining $P$ remains bounded if we remove  the half-space defined by
  $F$.
\end{definition}
See Figure~\ref{fig:ps01} for an  illustration of the applicability of
our definition.  The set   $\mathcal F:=\{F_1,F_2\}$ of facets  of the
polytope $P$ with the chosen simplex facet  $S$ in bounded position is
nonsimple, while the set $\mathcal F':=\{F_2,F_3\}$ is not.
\begin{remark}
  There  is always a projective   transformation  that puts a  simplex
  facet  $S$ of a $d$-polytope  $P$ into bounded  position, as long as
  there exists a facet of $P$ not adjacent to $S$, that is, if $P$
  is not a simplex.
  Without loss of generality we can therefore assume that a simplex facet
  is in bounded position, if the given polytope itself is not a simplex.
\end{remark}

In    the following  lemma we  show   that  we can   always apply  our
construction to a simplex facet in bounded position, regardless of the
choice of $\mathcal F$ and $\mathcal N$.
\begin{lemma}
  Let $P$ be a $d$-polytope and $S$ a  simplex facet of $P$ in bounded
  position.    Choose   two  disjoint   sets  $\mathcal   F,  \mathcal
  N\subseteq\adj{S}$.  Then the set $\nv{S}{\mathcal F,\mathcal N}{P}$
  is not empty.
\end{lemma}
\begin{proof}
  Let $F_1,\ldots, F_d$   be the facets   in $\adj S$   and number the
  vertices $p_1,\ldots,p_d$ of $S$ in such a way that $p_j$ is the
  unique vertex not lying in $F_j$ for $j=1,\ldots, d$.
  \par
  $S$ is a $(d-1)$-simplex and is  in bounded position with respect to $P$, so
  all hyperplanes in the  set $\hadj{\adj{F}}$ intersect  in a unique point
  $p\in H_S^-$.  Let $v_j := p_j-p$ for $j=1,\ldots, d$. Define $p'\in\R^d$ by
  \begin{align*}
    p' := p + \eps \left(\sum_{j : F_j \in \adj S \setminus ({\mathcal  F} \cup {\mathcal N})} v_j
      - \sum_{j : F_j \in {\mathcal N}} v_j \right) .
  \end{align*}
  For $\eps>0$ small   enough $p'$  is contained  in  $\nv{S}{\mathcal
    F,\mathcal N}{P}$.
\end{proof}

The following proposition tells which types of facets occur in 
$\pstack{S}{\mathcal F, \mathcal N}{P}$.
\begin{proposition} \label{prop:facets-ps} Let  $P$ be a $d$-polytope,
  $d\ge3$,  $S$  a simplex  facet  of   $P$ in  bounded position   and
  $\mathcal F,   \mathcal N\subseteq\adj S$,   $\mathcal F\cap\mathcal
  N=\emptyset$.   Assume  that $\mathcal F$   is nonsimple.   Then the
  following properties  of  $\pstack{S}{\mathcal F,   \mathcal  N}{P}$
  hold:
  \begin{compactenum}[\rm (1)]
  \item The relative interiors of $S$ and of the facets in $\mathcal N$ lie
    in the interior of $\pstack{S}{\mathcal F, \mathcal N}{P}$.
  \item Let $v$  be the   new  vertex added by the    pseudo-stacking.
    $\pstack{S}{\mathcal   F,  \mathcal N}{P}$     has facets of   the
    following four types:
    \begin{compactenum}[\rm (a)]
    \item Facets not   in $\{S\}\cup\mathcal N\cup\mathcal   F$ remain
      unchanged in $\pstack{S}{\mathcal F, \mathcal N}{P}$.
    \item  For  any ridge  $R$ between  $S$ and  a  facet  in  $\adj S
      \setminus ({\mathcal F} \cup {\mathcal N})$ we obtain  a new facet that
      is a pyramid over $R$ with apex $v$.
    \item For any ridge $R$ between a facet in $\mathcal N$ and one in
      $\facets(P)\setminus(\mathcal   F\cup\mathcal   N\cup\{S\})$  we
      obtain a new facet that is a pyramid over $R$ with apex $v$.
    \item For any facet $F\in\mathcal F$ we obtain a facet of the form
      $\pstack{S\cap   H_F}{\emptyset,\mathcal   N_F}{F}$    in    the
      pseudo-stacking, where we view $F$  as a  polytope in $H_F$  and
      define $\mathcal N_F:=\{N\cap F\mid N\in\mathcal N\}$.
    \end{compactenum}
  \end{compactenum}
\end{proposition}
\begin{proof}\hfill
  \begin{compactenum}
  \item $v$  lies beyond   $S$  and all   facets in  $\mathcal   N$ by
    definition. Hence, any ray from $v$ to a point $v'$ in the relative interior of one of
    these facets intersects $P$ in a segment with boundary points $v'$
    and some other point $v''$. The segment from $v$  to $v''$ is
    contained in $\conv(P,v)$. Hence, $v'$  must be in the interior of
    $\conv(P,v)$. This proves the first claim.
  \item The vertices of  $\conv(P,v)$ are a subset  of the vertices of
    $P$ and the vertex $v$. We show that all four facet types mentioned
    may occur in the pseudo-stacking, and no others.
    \begin{compactenum}
    \item The new vertex $v$ lies  beneath all facets not in $\mathcal
      F\cup\mathcal  N\cup\{S\}$. Hence the  facet hyperplanes of such
      facets remain valid  and facet defining  for the pseudo-stacking
      of $P$.   This completely describes  all  facets of $\conv(P,v)$
      that do not involve the vertex $v$.
    \item Let $F \not\in \mathcal F \cup \mathcal N$ be a facet adjacent to $S$
      and $R = F \cap S$. Then $v$ lies beneath $H_F$ and beyond $H_S$.
      Hence, any segment between $v$ and a point $v' \in R$ intersects $P$
      only in $v'$; on the other hand, a segment between $v$ and some
      $v'' \in F \setminus R$ must intersect the interior of $P$.
      In the convex hull of $v$ and $P$ we thus obtain a pyramid over $R$.
    \item As  above, there is a  segment from $v$  to any vertex  of a
      ridge  $R$ between a facet   $N\in\mathcal N$ and  a facet $F\in
      \facets(P)\setminus(\mathcal F\cup\mathcal N\cup\{S\})$, but not
      to any relative interior point of $F$. Hence, in the convex hull
      of $v$ and $P$ we again obtain a pyramid over $R$.
    \item Let  $F$ be a facet in   $\mathcal F$ and  $G$ the  ridge it
      shares  with $S$.  We look at  $F$ as a  polytope defined in the
      hyperplane  $H_F$.  The added  vertex    $v$ also lies in   that
      hyperplane. $G$ is a facet of $F$.  As $\mathcal F$ is nonsimple
      and $S$ a simplex facet, none of the facets of $F$ adjacent to
      $G$ is defined by a hyperplane coming from  a facet in $\mathcal
      F$.   Let $\mathcal N_F$  be the set of   ridges $F$ shares with
      some facet of $\mathcal N$. This is a set of facets of $F$. \\
      Hence, seen as a point in $H_F$, $v$ is beyond $G$, beyond all
      facets in $\mathcal N_F$ and beneath all other facet hyperplanes
      of $F$.  So $F$ is pseudo-stacked by $v$ above $G$ and
      $\mathcal N_F$.
    \end{compactenum}
    This completely describes all facets of $\conv(P,v)$ that involve the vertex $v$.
    Hence we   have  described    all   possible  facet    types    of
    $\conv(P,v)$.\qedhere
  \end{compactenum}  
\end{proof}
\begin{remark} \label{rem:bounded-pos-again}   Any simplex  facet   of
  $\pstack{S}{\mathcal F, \mathcal N}{P}$   of   type (b) or (c)    in
  Proposition~\ref{prop:facets-ps} is again in bounded position.
\end{remark}

The following theorem tells which $k$-faces of $P$ are also
$k$-faces of  $\pstack{S}{\mathcal  F,\mathcal N}{P}$, for $0\le  k\le
d-2$ (cf.\ also Gr\"un\-baum \cite[Sec.~5.2]{MR1976856}).
\begin{theorem} \label{thm:surviving-faces} Let $P$ be a $d$-polytope,
  $d\ge3$,   $S$ a simplex   facet  of  $P$  in  bounded  position and
  $\mathcal F,   \mathcal  N\subseteq\adj S$,  $\mathcal F\cap\mathcal
  N=\emptyset$.  Assume that $\mathcal F$ is nonsimple.
  \par
  A $k$-face $G$   of $P$, $0\le  k\le   d-2$ is again a   $k$-face of
  $\pstack{S}{\mathcal F,\mathcal N}{P}$ if   and only if there  is  a
  facet    $F\in\facets(P)  \setminus   (\mathcal F\cup\mathcal  N\cup
  \{S\})$ that contains $G$.
\end{theorem}
\begin{proof}
  If $G$ is contained in a facet
  $F \in \facets(P) \setminus (\mathcal F \cup \mathcal N \cup \{S\})$,
  then $F$ is also a facet of $\stdpstack$ by Proposition~\ref{prop:facets-ps},
  and $G$ is also a face of $\stdpstack$.
  \par
  So assume that $G$ is a face of $P$ all of whose incident facets are
  contained in $\mathcal A := \mathcal F \cup \mathcal N \cup \{S\}$.
  Then any ridge of $P$ containing $G$ is also
  only incident to facets in $\mathcal A$ and we can assume that $G$
  is a ridge of $P$. Because $S$ is a simplex, the intersection of $G$ and $S$
  contains a subridge. Hence, by nonsimplicity, there is at most one facet from
  $\mathcal F$ that contains $G$.
  \begin{compactenum}
  \item\label{van1} If there is no facet from $\mathcal F$ containing to
    $G$,  then $v$ lies beyond both facet hyperplanes defining $G$. Hence,  
    $G$ is in the interior  of $\stdpstack$.
  \item If one of the facets containing $G$ is $F \in \mathcal F$ and
    the   other is $S$,  then $F$   is  stacked above  $G$, hence  $G$
    vanishes in the interior of $F$.
  \item If one of the facets containing $G$ is $F \in \mathcal F$ and
    the other facet is   in  $\mathcal N$, then the    pseudo-stacking
    transforms $F$ into the facet $\pstack{S\cap F}{\emptyset,\mathcal N'}{F}$,
    where $\mathcal N' := \{ N \cap F \mid N \in \mathcal N \}$
    and $G$ is one of the facets in $\mathcal N'$ which vanish in
    the pseudo-stacking.\qedhere
  \end{compactenum}
\end{proof}
Restricting to the vertices of $P$ we have the following consequence.
\begin{corollary} \label{cor:vertices-plus-1}
  Let $P$ be a  $d$-polytope, $d\ge3$, $S$ a simplex  facet of $P$  in
  bounded position  and $\mathcal   F,  \mathcal  N\subseteq\adj   S$,
  $\mathcal F\cap\mathcal  N=\emptyset$.  Assume that $\mathcal  F$ is
  nonsimple and that any  vertex of $P$ is contained  in at  least one
  facet $F \in \facets(P) \setminus ({\mathcal F} \cup {\mathcal N} \cup \{S\})$.
  Then $f_0(\pstack{S}{\mathcal F, \mathcal N}{P})=f_0(P)+1$.
\end{corollary}
\begin{proof}
  By the previous  theorem all vertices  of  $P$ are also  vertices of
  $\stdpstack$. The pseudo-stacking operation adds one new vertex.
\end{proof}

The next theorems deal with the consequences of the pseudo-stacking 
operation on $k$-simpliciality and $h$-simplicity.
\begin{theorem} \label{thm:k-simplicial} Let  $P$ be a  $k$-simplicial
  $d$-polytope  for $d\ge3$ and $1\le k\le  d-2$  with a simplex facet
  $S$ in bounded position.  Let $\mathcal F, \mathcal N \subseteq \adj
  S$ be disjoint  sets.  Assume that $\mathcal  F$ is nonsimple.  Then
  $\pstack{S}{\mathcal F, \mathcal N}{P}$ is $k$-simplicial.
\end{theorem}
\begin{proof}
  Let $v$  be the  added  vertex in the  pseudo-stacking. We   use the
  characterisation  of the facets in Proposition~\ref{prop:facets-ps}.
  The  $k$-faces   in   facets  that  stay   unchanged   clearly  stay
  combinatorially equivalent to a $k$-simplex. The $k$-faces contained
  in facets that  are pyramids over ridges of  $P$  are either already
  faces  of $P$ or  pyramids over $(k-1)$-faces  of $P$ with apex $v$.
  Hence they are simplices.
  \par
  All remaining facets are obtained by pseudo-stacking a facet of $P$.
  However, only the case  $\mathcal  F=\emptyset$ occurs, that is,  in
  the pseudo-stacking of   a facet none  of its  facets gets  stacked.
  Therefore, the added facets are all pyramids over ridges, which preserve
  $k$-simpliciality for $k\le d-2$.
\end{proof}

Now we look more   closely at the types  of  edges that can  occur  in
$\pstack{S}{\mathcal F,\mathcal N}{P}$ and  determine their number. This 
is used to establish $2$-simplicity for the family of polytopes we 
construct in the next section. The two cases   
$\mathcal   N=\emptyset$   and   $\mathcal N\ne\emptyset$  are treated 
separately   in the   following  two propositions,  and
the latter case is further restricted to $|\mathcal N|=1$. 
\begin{definition}
  Let $P$ be a polytope.  For any face $F$  of $P$ let $\fdeg_P(F)$ be
  the number of facets of $P$ that contain $F$.
\end{definition}
\begin{proposition} \label{prop:new-edges}  Let $P$ be a $d$-polytope,
  $d\ge3$, $S$ a simplex facet of $P$ in bounded position, $\mathcal F
  \subseteq  \adj     S$  nonsimple   and    $v  \in   \nv{S}{\mathcal
    F,\emptyset}{P}$.
  \begin{compactenum}[\rm (1)]
  \item All edges of  $\pstack{S}{\mathcal  F}{P}$ are  of one of  the
    following two types:
    \begin{compactenum}[\rm (a)]
    \item edges $e$ in $P$ such that $e \subset F$ for a facet
      $F \in \facets(P) \setminus (\{S\} \cup \mathcal F)$;
    \item edges $e=[v,v']$ for every vertex $v' \in S$; in this case
      \[ \fdeg_{\pstack{S}{\mathcal F}{P}}(e) = \fdeg_S(v') . \]
    \end{compactenum}
  \item  $f_1(\pstack{S}{\mathcal F}{P}) =   f_1(P)  + d  -  f$, where
    $f=|\mathcal F|$ if $d=3$ and $f=0$ otherwise.
  \end{compactenum}
\end{proposition}
\begin{proof}\hfill
  \begin{compactenum}
  \item The edges of $\pstack{S}{\mathcal  F}{P}$ that are also  edges
    of  $P$     are  exactly    the   ones  described   in    (a),  by
    Theorem~\ref{thm:surviving-faces}.  On  the other hand, if an edge
    $e$ of $\pstack{S}{\mathcal F}{P}$ is not an edge of $P$, then one
    of its vertices is $v$  and the other vertex  $v'$ is contained in
    $S$, which is the only facet of $P$ whose vertices can be seen
    by $v$.  In this case,  every facet of $\pstack{S}{\mathcal F}{P}$
    that contains  $e$ corresponds to a ridge  $R \subset S$ such that
    $v' \in  R$ --- either the facet  is of type (b) or  of type (d) in
    Proposition~\ref{prop:facets-ps}.
  \item Since $S$  is a simplex  there are  $d$ edges of type  (b). If
    $d>3$ then, since $\mathcal F$ is nonsimple,  every edge of $P$ is
    contained  in at least  one facet not  in $\mathcal F \cup \{S\}$,
    hence there are $f_1(P)$  edges of type  (a). If $d=3$, then edges
    are  ridges and every   ridge $F \cap S$   for $F \in  \mathcal F$
    disappears -- namely, it is the base facet for the pseudo-stacking
    of $F$ (cf.~Proposition~\ref{prop:facets-ps}). In this case we
    have $f_1(P)-|\mathcal F|$ edges of type (a).\qedhere
  \end{compactenum}
\end{proof}

\begin{proposition}  \label{prop:new-edges-N}   Let    $P$  be       a
  $d$-polytope,  $d\ge3$,  $S$ a   simplex  facet of  $P$  in  bounded
  position, $\mathcal F  \subseteq  \adj S$ nonsimple, $\mathcal  N  =
  \{N\}$ with a facet $N \in  \adj S \setminus  \mathcal F$ and $v \in
  \nv{S}{\mathcal F,\mathcal N}{P}$.  Additionally, suppose that every
  vertex of $S$   and $N$ is contained   in  at least one   facet from
  $\facets(P) \setminus (\{S\} \cup \mathcal F \cup \mathcal N)$.
  \begin{compactenum}[\rm (1)]
  \item All edges of $\pstack{S}{\mathcal F,\mathcal N}{P}$ are of one
    of the following types:
    \begin{compactenum}[\rm (a)]
    \item edges $e$ in $P$ such that $e \subset  F$ for a facet $F \in
      \facets(P) \setminus (\{S\} \cup \mathcal F \cup \mathcal N)$;
    \item edges $e=[v,v']$ for every vertex $v'$ of either $S$ or $N$;
      in this case
      \[ \fdeg_{\pstack{S}{\mathcal F,\mathcal N}{P}}(e) = \fdeg_X(v')
      , \] where $X=S$ or $X=N$ respectively;
    \item edges $e=[v,v']$ for every vertex $v'$ of $S \cap N$.
    \end{compactenum}
  \item  $f_1(\pstack{S}{\mathcal F,  \mathcal N}{P})  = f_1(P)  + d +
    f_0(N) - f_0 (S \cap N) - f$,  where $f=|\mathcal F|$ if $d=3$ and
    $f=0$ otherwise.
  \end{compactenum}
\end{proposition}
\begin{proof}\hfill
  \begin{compactenum}
  \item  The description of types  of  edges is  complete  by the same
    argument as in the proof before; the  only difference is that here
    the vertex $v$ sees both facets $S$ and $N$. Also, if $v'$ is
    a vertex in either  $S$ or $N$,  the edge $[v,v']$ is contained in
    the facets  corresponding to the ridges  of  $S$, resp.\ $N$, that
    contain $v'$.
  \item Again, since $S$ is a simplex, $\pstack{S}{\mathcal F,\mathcal N}{P}$ has
    $d + f_0(N)  - f_0(S \cap N)$  edges of type (b)  or (c). If $d>3$
    then $\mathcal F$  being nonsimple ensures  that every edge of $P$
    is also an edge of $\pstack{S}{\mathcal F,\mathcal N}{P}$ of type (a). If
    $d=3$, again all edges between  facets in $\mathcal  F$ and $S$ or
    $N$ disappear; there  is  no  facet in  $\mathcal F$   that shares
    ridges with both $S$ and $N$, since if there were, we had a vertex
    contained only in facets from $\{S\} \cup \mathcal F \cup \mathcal
    N$.   Therefore exactly  $f$  edges from  $P$    are no  edges  in
    $\pstack{S}{\mathcal F,\mathcal N}{P}$.\qedhere
  \end{compactenum}
\end{proof}

We want  to show  in the next   section that  certain combinations  of
pseudo-stacking operations  preserve $2$-simplicity of  a $4$-polytope
$P$.  For this we have to count the number of facets a subridge of $P$
is  in.  We do  this with the  following  two propositions.  Again, we
consider  the  two    cases  $\mathcal  N=\emptyset$    and  $\mathcal
N\ne\emptyset$ separately and restrict the latter to $|\mathcal N|=1$.
\begin{proposition}   \label{prop:subridge-deg}   Let  $P$     be    a
  $d$-polytope, $d\ge3$, $S$   a  simplex  facet  of $P$  in   bounded
  position, $\mathcal F \subseteq \adj S$ nonsimple and $G$ a subridge
  of $P$ with $G \subset S$.  Define $\varphi :=  |\{ F \in \mathcal F
  \, | \, G \subset F \}|$.
  \par
  Then   $0 \leq \varphi      \leq  2$  and    $G$  is  a  face     of
  $\pstack{S}{\mathcal         F}{P}$;     furthermore,             \[
  \fdeg_{\pstack{S}{\mathcal F}{P}}(G) = \fdeg_P(G) + 1 - \varphi . \]
\end{proposition}
\begin{proof}
  $G$ is a ridge of $S$,  hence there are  two facets $G_1,G_2$ of $S$
  such that $G_1 \cap G_2 = G$.  Since $G_1,G_2$ are ridges of $P$, if
  $G$ is contained  in some facet $F \in  \mathcal F$, then $F  \cap S
  \in \{ G_1,G_2 \}$. Therefore at most $2$ facets in $\mathcal F$ can
  possibly contain $G$.
  \par
  If $G$  was   not a face   of $\pstack{S}{\mathcal  F}{P}$ then   by
  Theorem~\ref{thm:surviving-faces}   it   would be  contained  in two
  facets $F_1,F_2 \in \mathcal F$.  Since $S$ is  a simplex, $F_1 \cap
  F_2$ was  a ridge containing $G$,  in contradiction  to $\mathcal F$
  being nonsimple.
  \par
  Now count  the   number of  facets  of $\pstack{S}{\mathcal   F}{P}$
  containing $G$.   All facets of $P$ not  in $\mathcal F  \cup \{S\}$
  remain     facets          of     $\pstack{S}{\mathcal  F}{P}$    by
  Proposition~\ref{prop:facets-ps};  therefore, $\fdeg_P(G)-\varphi-1$
  facets still contain $G$  in $\pstack{S}{\mathcal F}{P}$. All facets
  in   $\mathcal F$    containing     $G$   also  stay    facets    of
  $\pstack{S}{\mathcal F}{P}$ containing  $G$. Additionally, for every
  one of the  two ridges $R$  of $P$ with $G  \subset R \subset S$, we
  either pseudo-stacked the facet $F$ with $F \cap S = R$ (if it is in
  $\mathcal F$) or we get a new facet, which is a pyramid over $R$; in
  the latter case, we have not counted it yet.
  \par
  In  total,      there are  $\fdeg_P(G)-\varphi-1+\varphi+2-\varphi =
  \fdeg_P(G)+1-\varphi$  facets  of $\pstack{S}{\mathcal   F}{P}$ that
  contain $G$.
\end{proof}

\begin{proposition}    \label{prop:subridge-deg-N}   Let  $P$,    $S$,
  $\mathcal      F$,          $G$      and      $\varphi$   as      in
  Proposition~\ref{prop:subridge-deg}.  Choose  a facet  $N \in \adj S
  \setminus \mathcal F$ and set $\mathcal N = \{N\}$.
  \begin{compactenum}[\rm (1)]
  \item If $G \subset S$ and $G \not\subset N$, then $G$  is a face of
    $\pstack{S}{\mathcal  F,\mathcal  N}{P}$             and        \[
    \fdeg_{\pstack{S}{\mathcal F,\mathcal N}{P}}(G) = \fdeg_P(G) + 1 -
    \varphi . \]
  \item If $G \subset N$ and every facet containing $G$ is either
    $S$, $N$ or   in   $\mathcal F$,  then $G$  is   not a   face   of
    $\pstack{S}{\mathcal F,\mathcal N}{P}$;  otherwise, $G$ is  a face
    of $\pstack{S}{\mathcal F,\mathcal N}{P}$ and
    \[  \fdeg_{\pstack{S}{\mathcal F,\mathcal N}{P}}(G) = \fdeg_P(G) +
    1 - \varepsilon - \varphi , \] where $\varepsilon=1$ if $G \subset
    S$ and $\varepsilon=0$ if not.
  \end{compactenum}
\end{proposition}
\begin{proof}
  (1)     can        be    shown   in      the       same     way   as
  Proposition~\ref{prop:subridge-deg}, since   the   facet $N$ has  no
  influence whatsoever on $G$.
  \par
  Suppose      $G$       is      contained        in    $N$.        By
  Theorem~\ref{thm:surviving-faces},       $G$  is   a    face      of
  $\pstack{S}{\mathcal F,\mathcal N}{P}$ if and only if there is facet
  containing $G$ that is not in $\mathcal F \cup \{ S,N \}$. In this
  case, the number of facets of $\pstack{S}{\mathcal F,\mathcal N}{P}$
  can be  computed in the same  way  as in the  previous proof, except
  that  if $G \subset  S$, an additional  facet containing $G$ (namely
  $S$) disappears.
\end{proof}

\begin{remark}\label{rem:generalised-pseudo}
  The construction given in this section can be further generalised:
  \begin{compactenum}
  \item Nonsimplicity of the set $\mathcal F$ is not necessary.
  \item The facet $S$ need not be a simplex.
  \item In this case the construction can be  adapted to only apply to
    a part of $S$.
  \end{compactenum}
  Several  new families of  polytopes arise via these generalisations.
  However,  keeping  track of  the changes  in  the $f$-vector is more
  subtle in these cases.  We  do not need them here,  so we omit their
  treatment.
\end{remark}

%%%%%%%%%%%%%%%%%%%%%%%%%%%%%%%%%%%%%%%%%%%%%%%%%%%%%%%%%%%%%%%%%%%%%%%%%%%%%%%

\section{Generating $2$-simple and $2$-simplicial $4$-polytopes} \label{sec:constr-2s2s}

In  this  part we  use  the  pseudo-stacking operation defined  in the
previous  section   to   extend elementary  $2$-simple  and    
$2$-simplicial $4$-polytopes  while   maintaining    these
properties.    As  we  have   seen in  Theorem~\ref{thm:k-simplicial},
pseudo-stacking preserves $2$-simpliciality.   The more difficult part
is to also preserve $2$-simplicity; we accomplish this by applying the
operation five, resp.\  four times in a  suitable way, such that edges
whose degree initially  increases  are  again  contained in only   $3$
facets at the end.

To simplify the task of keeping track of  the facets generated in each
step of the   construction  we introduce  two ways  of  distinguishing
certain facets.
\begin{definition}
  Let $P$ be a $d$-polytope and $S$ a facet of $P$.
  \begin{compactenum}[\rm (1)]
  \item If the vertices $v_0,\ldots,v_k$ all lie
    in one facet of $P$ and define this facet uniquely, then we denote
    this facet by $F(v_0,\ldots,v_k)$.
  \item If $e_0,\ldots,e_m$ are edges  of a facet $G$ of $S$, 
    then we denote the unique facet of $P$ adjacent to $S$ via 
    $G$ by $F_S(e_0,\ldots,e_m)$.
  \end{compactenum}
\end{definition}

%%%%%%%%%%%%%%%%%%%%%%%%%%%%%%%%%%%%%%%%%%%%%%%%%%

\subsection{The first construction}
\begin{figure}[p]
  \centering
  \psfrag{e0}[tr][tr]{$e_0$}
  \psfrag{e1}[tl][tl]{$e_1$}
  \psfrag{e2}[tr][tr]{$e_2$}
  \psfrag{e3}[br][br]{$e_3$}
  \psfrag{e4}[tr][tr]{$e_4$}
  \psfrag{e5}[bl][bl]{$e_5$}
  \psfrag{e6}[br][br]{$e_6$}
  \psfrag{e7}[tr][tr]{$e_7$}
  \psfrag{e8}[bl][bl]{$e_8$}
  \psfrag{e9}[br][br]{$e_9$}
  \psfrag{v0}[t][t]{$v_0$}
  \psfrag{v1}[l][l]{$v_1$}
  \psfrag{v2}[bl][bl]{$v_2$}
  \psfrag{v3}[bl][bl]{$v_3$}
  \psfrag{v4}[br][br]{$v_4$}
  \subfigure[Labelling of the vertices and edges of the simplex $S$.]{
    \begin{minipage}{.47\textwidth}\centering
      \includegraphics[width=.85\textwidth]{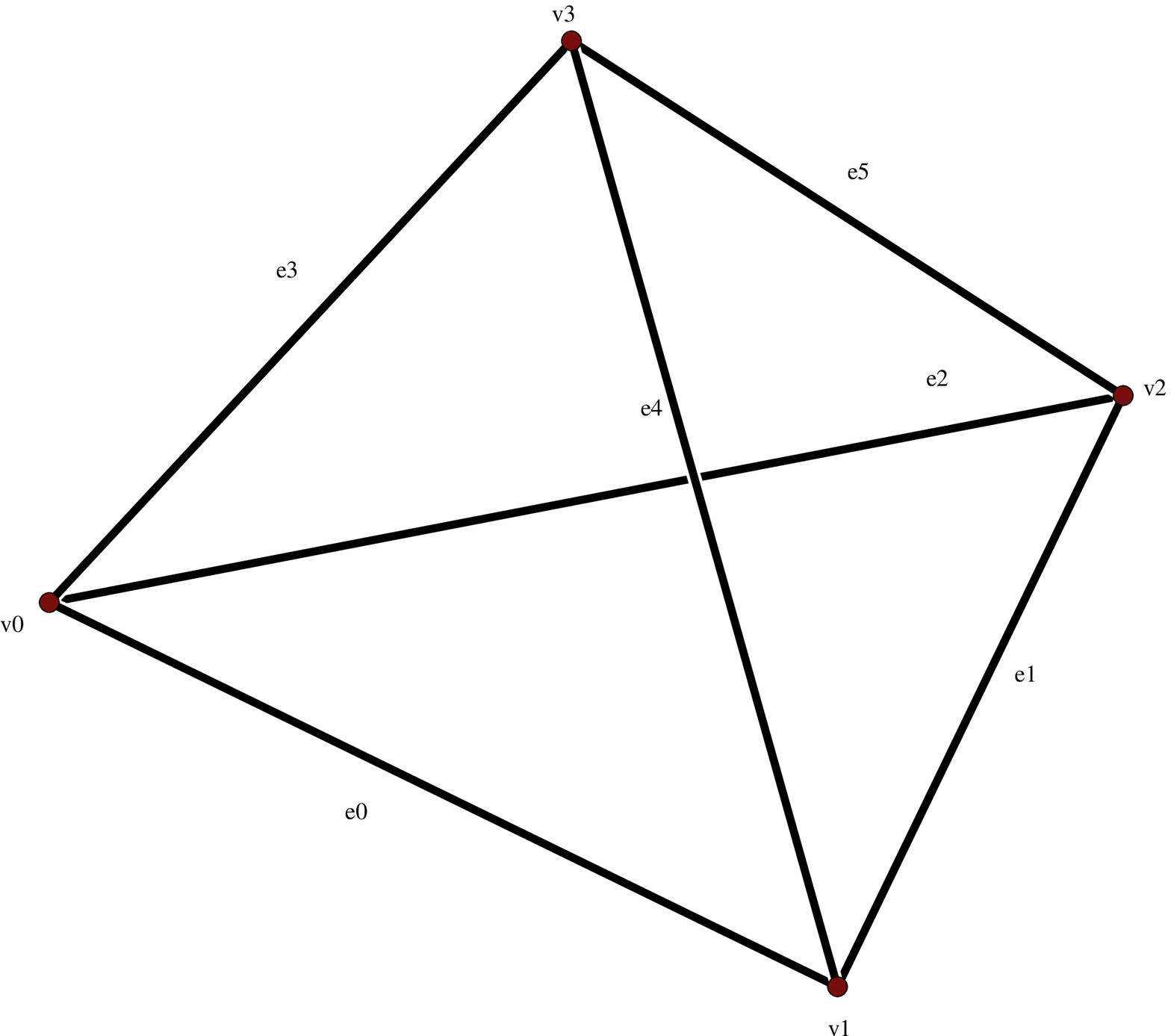}
      \label{fig:simplex-edge-labelling}
    \end{minipage}}
  \hfill
  \subfigure[Part of the Schlegel diagram of $P^{(1)}$;
  the new edges $e_6,\ldots,e_9$ (partly labelled) are drawn as dashed lines.]{
    \begin{minipage}{.47\textwidth}\centering
      \includegraphics[width=.85\textwidth]{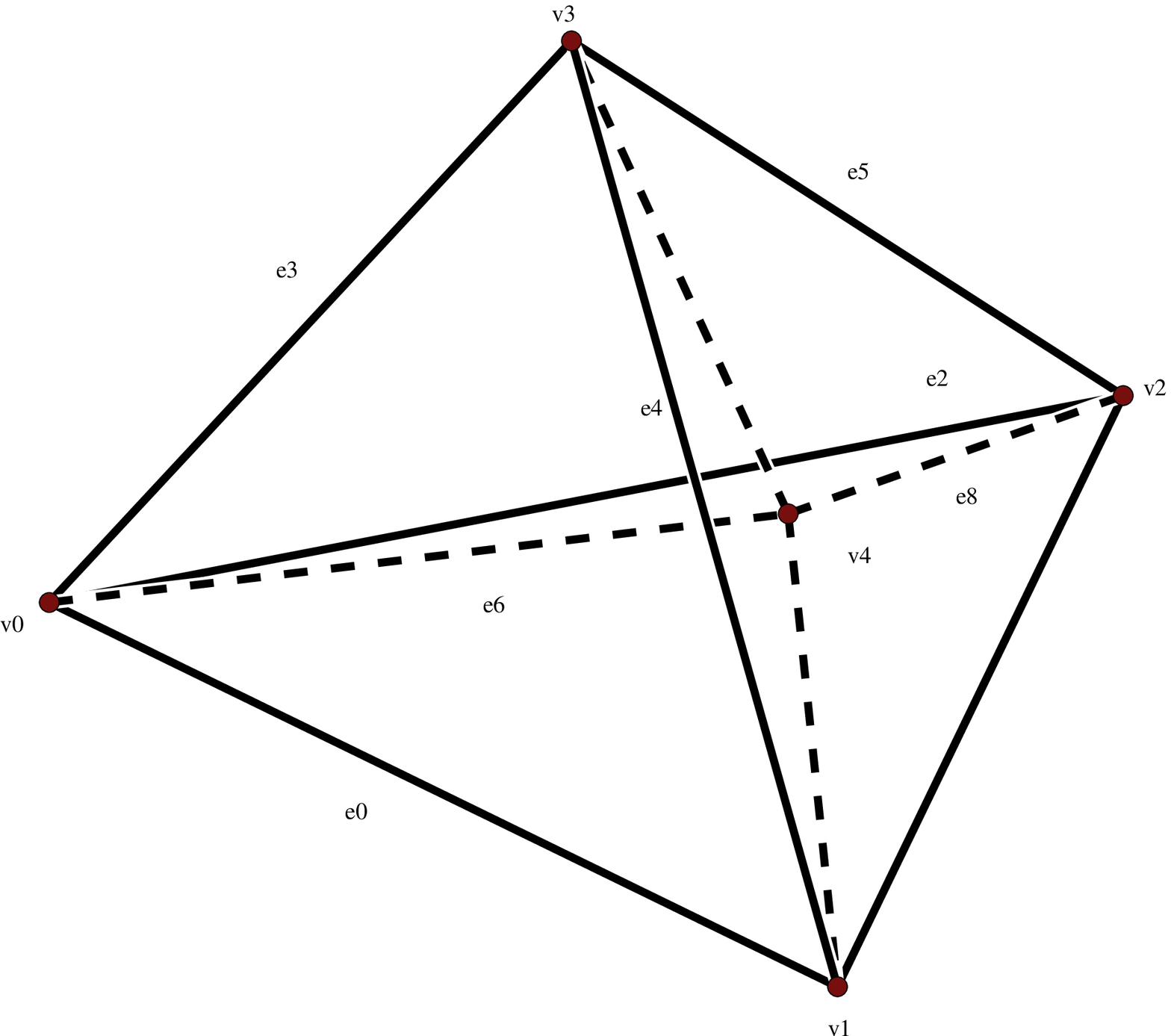}
      \label{fig:first-constr-1}
    \end{minipage}
  }
  
  \centering
  \psfrag{e0}[tr][tr]{$e_0$}
  \psfrag{e1}[tl][tl]{$e_1$}
  \psfrag{e2}[tr][tr]{$e_2$}
  \psfrag{e3}[br][br]{$e_3$}
  \psfrag{e4}[tr][tr]{$e_4$}
  \psfrag{e5}[bl][bl]{$e_5$}
  \psfrag{e6}[br][br]{$e_6$}
  \psfrag{e7}[tr][tr]{$e_7$}
  \psfrag{e8}[bl][bl]{$e_8$}
  \psfrag{e9}[br][br]{$e_9$}
  \psfrag{v0}[t][t]{$v_0$}
  \psfrag{v1}[l][l]{$v_1$}
  \psfrag{v2}[bl][bl]{$v_2$}
  \psfrag{v3}[bl][bl]{$v_3$}
  \psfrag{v4}[br][br]{$v_4$}
  \psfrag{v5}[br][br]{$v_5$}
  \subfigure[Part of the Schlegel diagram of $P^{(2)}$;
  note that $v_1\ldots,v_5$ define a bipyramid over a triangle,
  which occurred by stacking the facet $F(v_1,v_2,v_3,v_4)$ of $P^{(1)}$.]{
    \begin{minipage}{.47\textwidth}\centering
      \includegraphics[width=.85\textwidth]{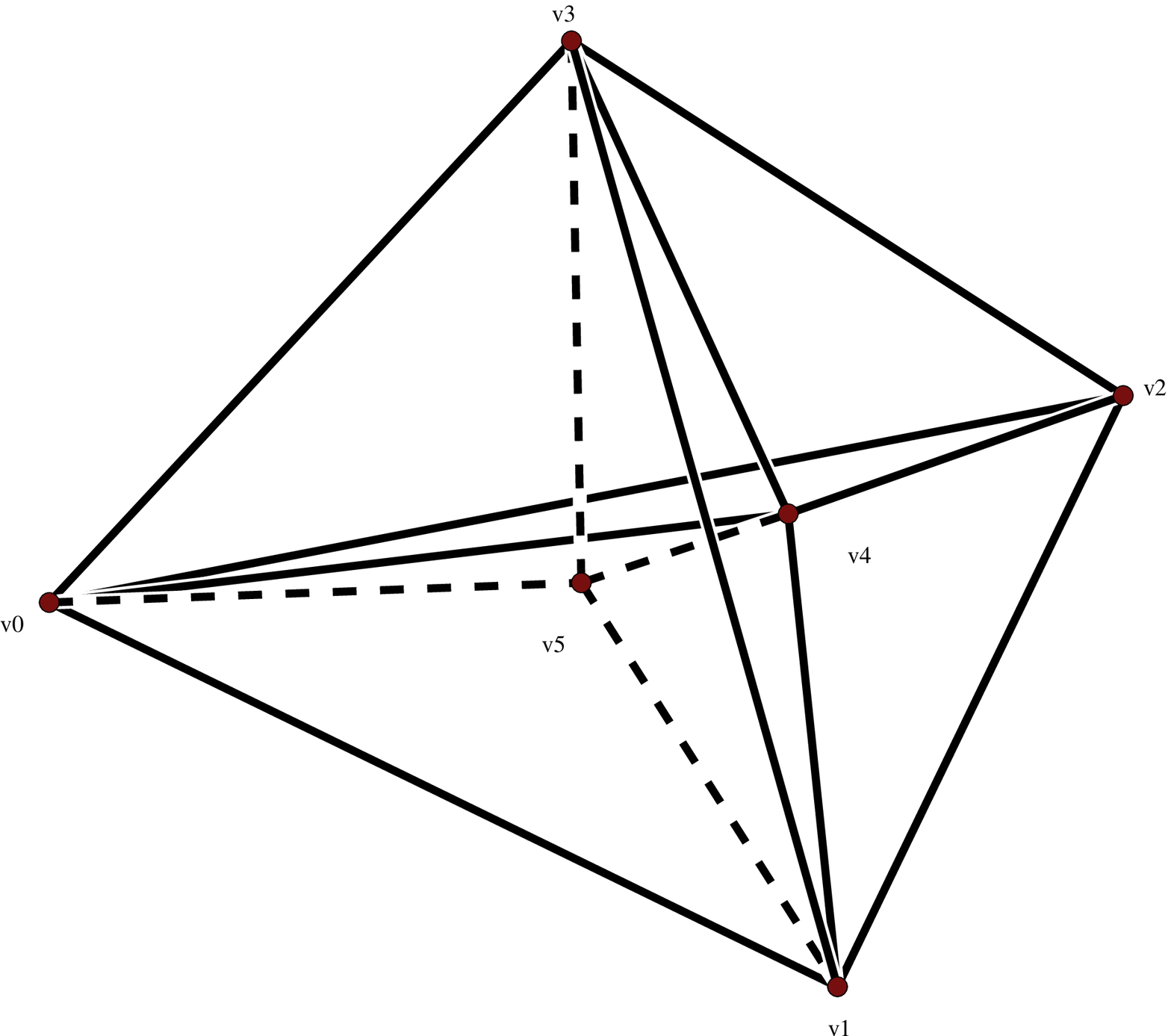}
      \label{fig:first-constr-2}
    \end{minipage}}
  \hfill
  \psfrag{e16}[tl][tl]{$e_{15}$}
  \psfrag{e17}[tr][tr]{$e_{16}$}
  \psfrag{e3}[br][br]{$e_3$}
  \psfrag{e4}[tr][tr]{$e_4$}
  \psfrag{e5}[bl][bl]{$e_5$}
  \psfrag{e6}[br][br]{$e_6$}
  \psfrag{e7}[tr][tr]{$e_7$}
  \psfrag{e8}[bl][bl]{$e_8$}
  \psfrag{e9}[br][br]{$e_9$}
  \psfrag{v0}[t][t]{$v_0$}
  \psfrag{v1}[l][l]{$v_1$}
  \psfrag{v2}[bl][bl]{$v_2$}
  \psfrag{v3}[bl][bl]{$v_3$}
  \psfrag{v4}[br][br]{$v_4$}
  \psfrag{v5}[br][br]{$v_5$}
  \psfrag{v6}[br][br]{$v_6$}
  \subfigure[Part of the Schlegel diagram of $P^{(3)}$;
  the shaded facets are the base facets
  for the last two pseudo-stacking steps.]{
    \begin{minipage}{.47\textwidth}\centering
      \includegraphics[width=.85\textwidth]{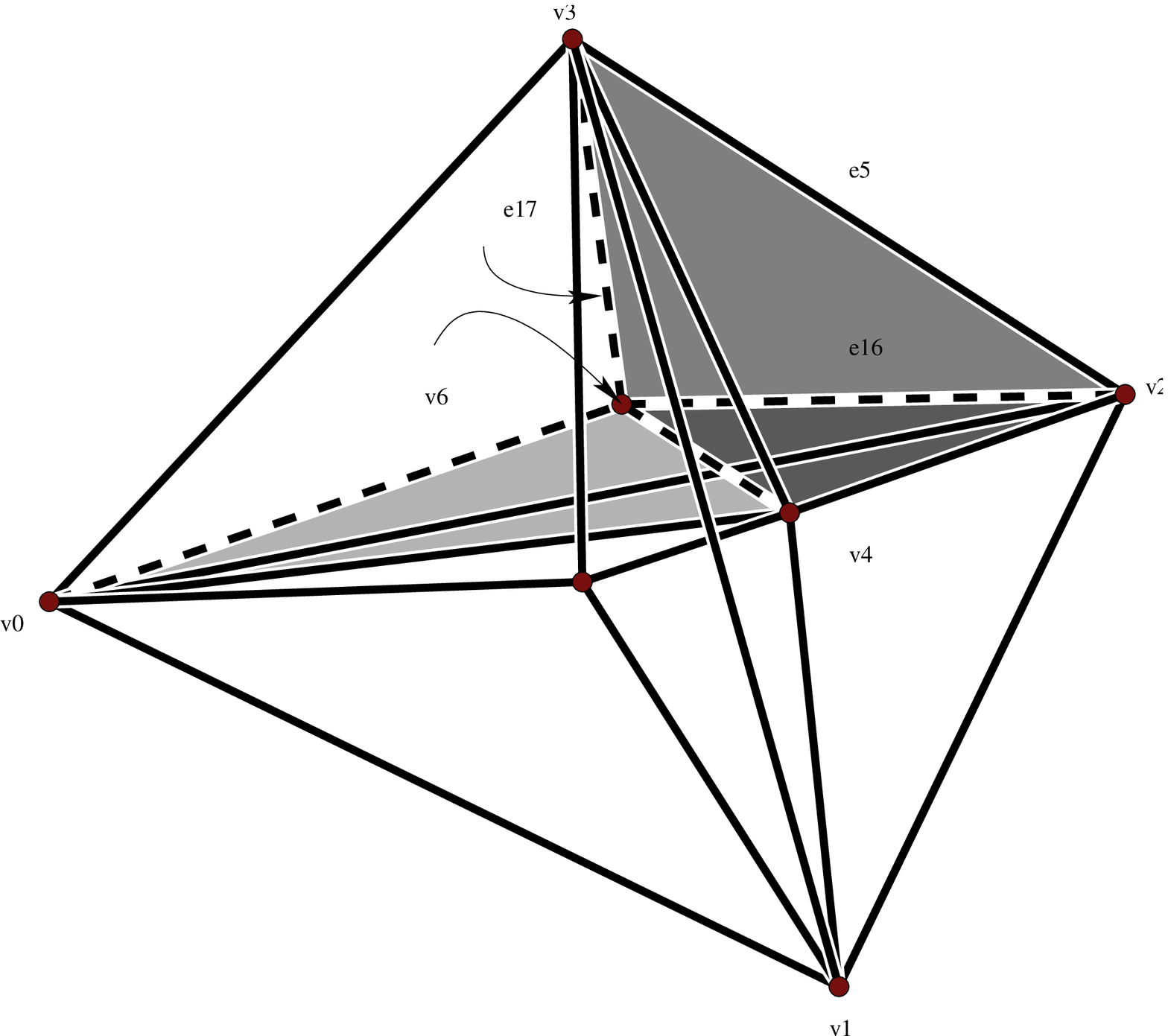}
      \label{fig:first-constr-3}
    \end{minipage}
  }
  
  \psfrag{e16}[tl][tl]{$e_{15}$}
  \psfrag{e17}[tr][tr]{$e_{16}$}
  \psfrag{e3}[br][br]{$e_3$}
  \psfrag{e4}[tr][tr]{$e_4$}
  \psfrag{e5}[bl][bl]{$e_5$}
  \psfrag{e6}[br][br]{$e_6$}
  \psfrag{e7}[tr][tr]{$e_7$}
  \psfrag{e8}[bl][bl]{$e_8$}
  \psfrag{e9}[br][br]{$e_9$}
  \psfrag{v0}[t][t]{$v_0$}
  \psfrag{v1}[l][l]{$v_1$}
  \psfrag{v2}[bl][bl]{$v_2$}
  \psfrag{v3}[bl][bl]{$v_3$}
  \psfrag{v4}[br][br]{$v_4$}
  \psfrag{v5}[br][br]{$v_5$}
  \psfrag{v6}[br][br]{$v_6$}
  \psfrag{v7}[br][br]{$v_7$}
  \subfigure[Part of the Schlegel diagram of $P^{(4)}$
  showing the pseudo-stacking of the facet $F(v_2,v_3,v_4,v_6)$ of $P^{(3)}$.]{
    \begin{minipage}{.47\textwidth}\centering
      \includegraphics[width=.85\textwidth]{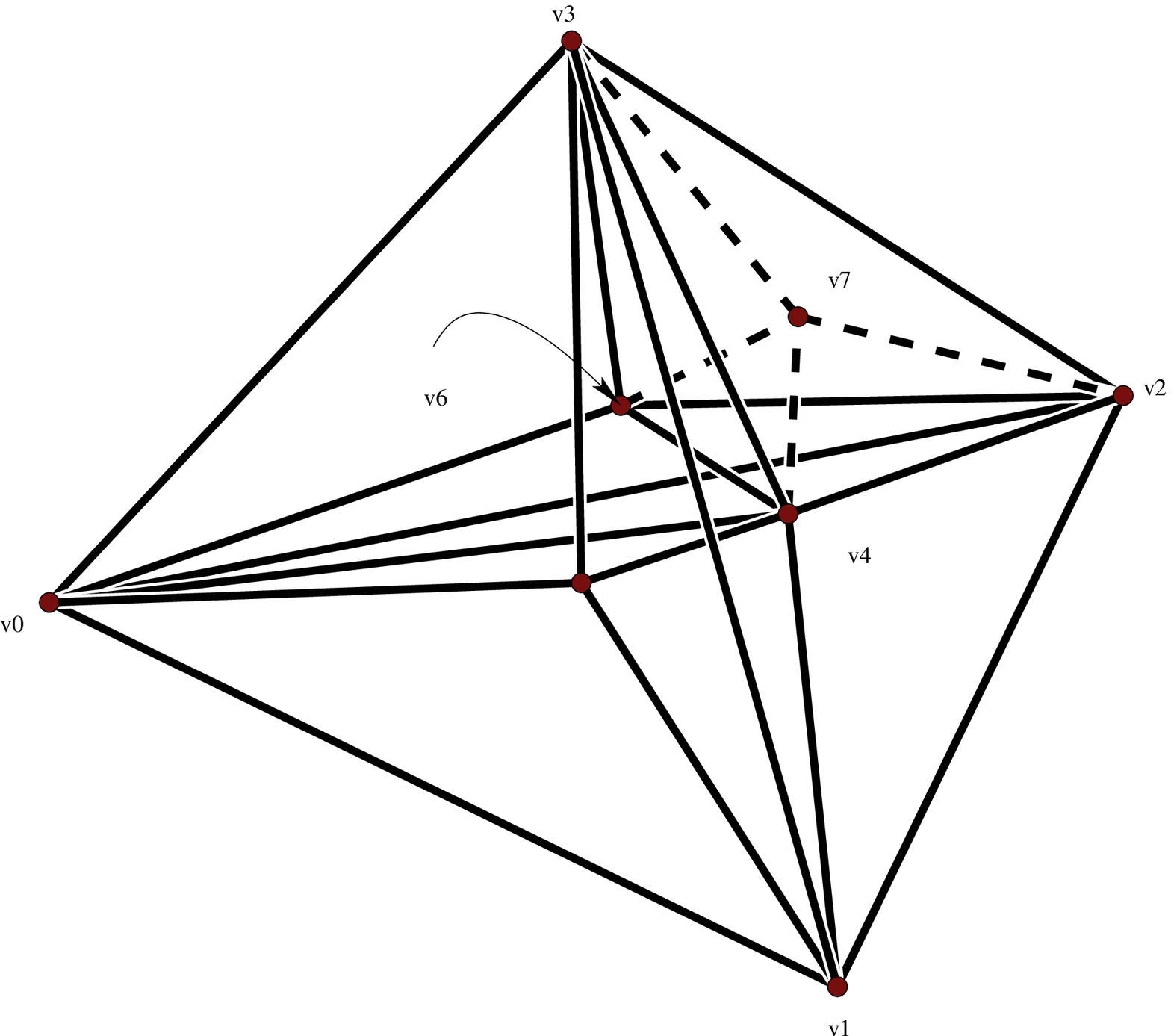}
      \label{fig:first-constr-4}
    \end{minipage}
  }
  \hfill
  \psfrag{v0}[t][t]{$v_0$}
  \psfrag{v1}[l][l]{$v_1$}
  \psfrag{v2}[bl][bl]{$v_2$}
  \psfrag{v3}[bl][bl]{$v_3$}
  \psfrag{v4}[br][br]{$v_4$}
  \psfrag{v5}[br][br]{$v_5$}
  \psfrag{v6}[br][br]{$v_6$}
  \psfrag{v7}[br][br]{$v_7$}
  \psfrag{v8}[br][br]{$v_8$}
  \subfigure[Part of the Schlegel diagram of $P^{(5)}$
  showing the pseudo-stacking of the facet $F(v_0,v_2,v_4,v_6)$ of $P^{(4)}$.]{
    \begin{minipage}{.47\textwidth}
      \centering
      \includegraphics[width=.85\textwidth]{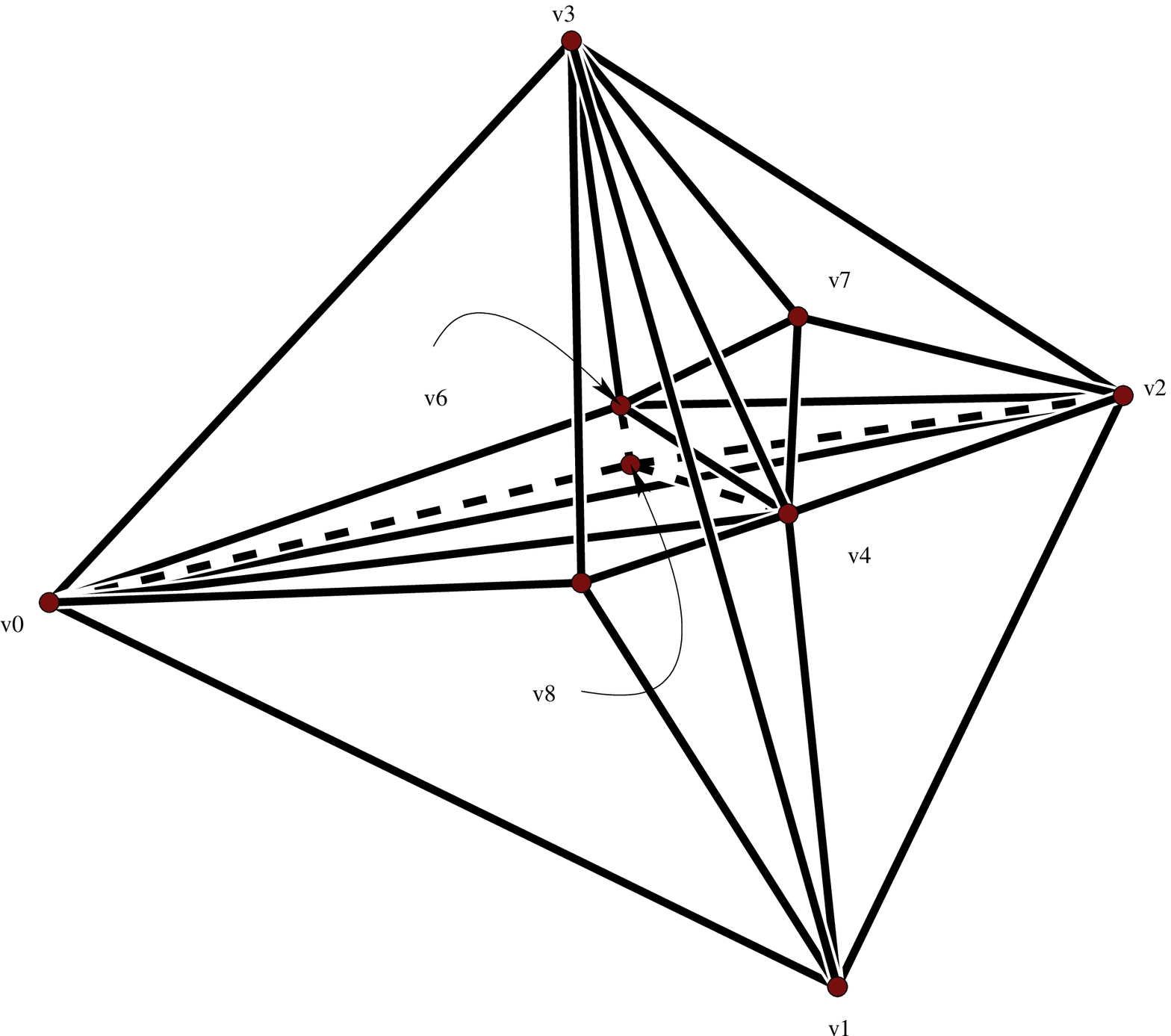}
      \label{fig:first-constr-5}
    \end{minipage}
  }

  \caption{The first construction.}
  \label{fig:the-first-construction}
\end{figure}

Let  $P$ be a $4$-polytope and  $S$ a simplex facet  of $P$ in bounded
position.  $S$ has   six  edges and  four  vertices,  which  we  label
$e_0,e_1,\ldots,e_5$    and    $v_0, v_1,    v_2,  v_3$   according to
Figure~\ref{fig:simplex-edge-labelling}.

We proceed with the five pseudo-stacking steps described in Table~\ref{tab:first}.
In Step $i$ we construct the polytope $P^{(i)} := \pstack{S_i}{{\mathcal F}_i}{P^{(i-1)}}$,
with $P^{(0)} := P$ and the parameters as given
in  Table~\ref{tab:first}.   For  every step  we  additionally give  a
complete       list       of    facets     of      type     (b)     in
Proposition~\ref{prop:facets-ps}, as well as of edges that either show
up  by  Proposition~\ref{prop:new-edges} or   change  their  degree by
Proposition~\ref{prop:subridge-deg}; note that the  only edges of  $P$
that change   degree during the process are   $e_3$, $e_4$  and $e_5$,
whose  degrees  in  $P$  are  denoted   by  $d_3$,   $d_4$  and  $d_5$
respectively.

Due to Remark~\ref{rem:bounded-pos-again} the base facets in
each step are in bounded   position; additionally, it can be  verified
easily that the sets ${\mathcal F}_i$ are nonsimple.
\begin{table}[tb]
  \footnotesize
  \caption{Summary of parameters and involved faces for the first construction.}
  \label{tab:first}
    \begin{tabular}{cccccc}
      step $i$ & base facet $S_i$ & facets in ${\mathcal F}_i$ & new
      facets & changing edges & degree \\ \hline
      $1$                             & $S$                    & $F_S(e_0,e_1,e_2)$           & $F(v_0,v_1,v_3,v_4)$ & $e_6:=[v_0,v_4]$    & $3$    \\
      (Fig.~\ref{fig:first-constr-1}) &                        &                              & $F(v_1,v_2,v_3,v_4)$ & $e_7:=[v_1,v_4]$    & $3$    \\
      \mbox{}                         &                        &                              & $F(v_0,v_2,v_3,v_4)$ & $e_8:=[v_2,v_4]$    & $3$    \\
      \mbox{}                         &                        &                              &                      & $e_9:=[v_2,v_4]$    & $3$    \\
      \mbox{}                         &                        &                              &                      & $e_3$               & $d_3+1$    \\
      \mbox{}                         &                        &                              &                      & $e_4$               & $d_4+1$    \\
      \mbox{} & & & & $e_5$ & $d_5+1$ \\ \hline
      $2$                             & $F(v_0,v_1,v_3,v_4)$   & $F(v_1,v_2,v_3,v_4)$         & $F(v_0,v_3,v_4,v_5)$ & $[v_0,v_5]$         & $3$    \\
      (Fig.~\ref{fig:first-constr-2}) &                        & $F_{S_2}(e_0,e_3,e_4)$       & $F(v_0,v_1,v_4,v_5)$ & $[v_1,v_5]$         & $3$    \\
      \mbox{}                         &                        &                              &                      & $[v_3,v_5]$         & $3$    \\
      \mbox{}                         &                        &                              &                      & $[v_4,v_5]$         & $3$    \\
      \mbox{}                         &                        &                              &                      & $e_4$               & $d_4$    \\
      \mbox{} & & & & $e_6$ & $4$ \\ \hline
      $3$                             & $F(v_0,v_2,v_3,v_4)$   & $F(v_0,v_3,v_4,v_5)$         & $F(v_2,v_3,v_4,v_6)$ & $e_{14}:=[v_0,v_6]$ & $3$    \\
      (Fig.~\ref{fig:first-constr-3}) &                        & $F_{S_3}(e_2,e_3,e_5)$       & $F(v_0,v_2,v_4,v_6)$ & $e_{15}:=[v_2,v_6]$ & $3$    \\
      \mbox{}                         &                        &                              &                      & $e_{16}:=[v_3,v_6]$ & $3$    \\
      \mbox{}                         &                        &                              &                      & $e_{17}:=[v_4,v_6]$ & $3$    \\
      \mbox{}                         &                        &                              &                      & $e_3$               & $d_3$    \\
      \mbox{} & & & & $e_8$ & $4$ \\ \hline
      $4$                             & $F(v_2,v_3,v_4,v_6)$   & $F(v_1,v_2,v_3,v_4)$         & $F(v_3,v_4,v_6,v_7)$ & $[v_2,v_7]$         & $3$    \\
      (Fig.~\ref{fig:first-constr-4}) &                        & $F_{S_4}(e_5,e_{16},e_{17})$ & $F(v_2,v_4,v_6,v_7)$ & $[v_3,v_7]$         & $3$    \\
      \mbox{}                         &                        &                              &                      & $[v_4,v_7]$         & $3$    \\
      \mbox{}                         &                        &                              &                      & $[v_6,v_7]$         & $3$    \\
      \mbox{}                         &                        &                              &                      & $e_5$               & $d_5$    \\
      \mbox{} & & & & $e_{17}$ & $4$ \\ \hline
      $5$                             & $F(v_0,v_2,v_4,v_6)$   & $F(v_0,v_3,v_4,v_5,v_6)$     & $F(v_0,v_2,v_6,v_8)$ & $[v_0,v_8]$         & $3$    \\
      (Fig.~\ref{fig:first-constr-5}) &                        & $F(v_2,v_4,v_6,v_7)$         &                      & $[v_2,v_8]$         & $3$    \\
      \mbox{}                         &                        & $F_{S_5}(e_2,e_6,e_8)$       &                      & $[v_4,v_8]$         & $3$    \\
      \mbox{}                         &                        &                              &                      & $[v_6,v_8]$         & $3$    \\
      \mbox{}                         &                        &                              &                      & $e_6$               & $3$    \\
      \mbox{}                         &                        &                              &                      & $e_8$               & $3$    \\
      \mbox{} & & & & $e_{17}$ & $3$ \\ \hline
    \end{tabular}
\end{table}

\begin{definition}
  Let $P$ be   a $4$-polytope with   a  simplex facet $S$ in   bounded
  position.  We denote the polytope $P^{(5)}$ obtained by applying the
  five steps of Table~\ref{tab:first} by $\first(P;S)$. 
\end{definition}
\begin{remark}\label{rem:ambiguous}
  Note that $\first(P;S)$ implicitly depends on a labelling of the vertices 
  of $S$. Choosing different labellings may result in combinatorially 
  different polytopes for the same choice of $P$ and $S$.
\end{remark}

\begin{lemma}\label{lemma:bounded-1}
  Let $P$  be   a $4$-polytope with  a  simplex  facet $S$ in  bounded
  position and the vertices of $S$  numbered in arbitrary order.  Then
  $\first(P;S)$   has   again  simplex facets    in  bounded position.
  Additionally,  all edges of $P$ are  still present in $\first(P;S)$,
  their   degrees  remain unchanged  and all  new   edges  have degree
  $3$.\qed
\end{lemma}

%%%%%%%%%%%%%%%%%%%%%%%%%%%%%%%%%%%%%%%%%%%%%%%%%%

\subsection{A second construction}

Let  $Q$  be  a $4$-polytope  with  a  simplex  facet   $S$ in bounded
position.  We label the vertices and edges  of $S$ in  the same way as
in the previous section (see Figure~\ref{fig:simplex-edge-labelling}).

The construction is  described by Table~\ref{tab:second},  in the same
way as  before.   The main difference is   that in Step  $3$  the set
${\mathcal    N}_3$  is   not empty,   so    in general  we  construct
$Q^{(i)}:=\pstack{S_i}{{\mathcal F}_i,{\mathcal N}_i}{Q^{(i-1)}}$ with
the given parameters. Note that one effect of this is that in Step $3$
one  edge disappears  by Theorem~\ref{thm:surviving-faces};  also, the
only edges  of $Q$ changing their  degree are  $e_1$, $e_2$ and $e_5$,
whose degrees are labelled accordingly.
\begin{figure}[bt]
  \subfigure[Part of the Schlegel diagram of $Q^{(1)}$
      with the labelling inherited from the facet $S$ of $Q$.]{
    \psfrag{e0}[tr][tr]{$e_0$}
    \psfrag{e1}[tl][tl]{$e_1$}
    \psfrag{e2}[tr][tr]{$e_2$}
    \psfrag{e3}[br][br]{$e_3$}
    \psfrag{e4}[tr][tr]{$e_4$}
    \psfrag{e5}[bl][bl]{$e_5$}
    \psfrag{e6}[br][br]{$e_6$}
    \psfrag{e7}[tr][tr]{$e_7$}
    \psfrag{e8}[bl][bl]{$e_8$}
    \psfrag{e9}[br][br]{$e_9$}
    \psfrag{v0}[t][t]{$v_0$}
    \psfrag{v1}[l][l]{$v_1$}
    \psfrag{v2}[bl][bl]{$v_2$}
    \psfrag{v3}[bl][bl]{$v_3$}
    \psfrag{v4}[bl][bl]{$v_4$}
    \includegraphics[width=.40\textwidth]{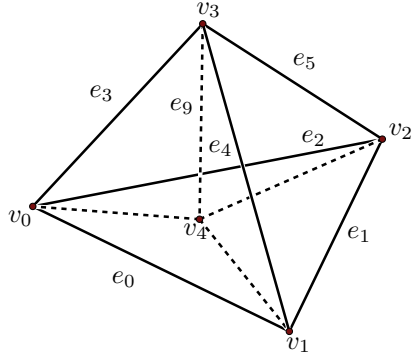}
     \label{fig:second-constr-1}}
\hspace{.05\textwidth}
  \subfigure[Part of the Schlegel diagram of $Q^{(2)}$,
      where the dashed edges were obtained by pseudo-stacking
      the facet $F(v_1,v_2,v_3,v_4)$ from the previous step.]{
    \psfrag{e0}[tr][tr]{$e_0$}
    \psfrag{e1}[tl][tl]{$e_1$}
    \psfrag{e2}[tr][tr]{$e_2$}
    \psfrag{e3}[br][br]{$e_3$}
    \psfrag{e4}[tr][tr]{$e_4$}
    \psfrag{e5}[bl][bl]{$e_5$}
    \psfrag{e6}[br][br]{$e_6$}
    \psfrag{e7}[tr][tr]{$e_7$}
    \psfrag{e8}[bl][bl]{$e_8$}
    \psfrag{e9}[br][br]{$e_9$}
    \psfrag{e10}[tr][tr]{$e_{10}$}
    \psfrag{e11}[tr][tr]{$e_{11}$}
    \psfrag{e12}[tr][tr]{$e_{12}$}
    \psfrag{e13}[tr][tr]{$e_{13}$}
    \psfrag{v0}[t][t]{$v_0$}
    \psfrag{v1}[l][l]{$v_1$}
    \psfrag{v2}[bl][bl]{$v_2$}
    \psfrag{v3}[bl][bl]{$v_3$}
    \psfrag{v4}[bl][bl]{$v_4$}
    \psfrag{v5}[bl][bl]{$v_5$}
    \includegraphics[width=.40\textwidth]{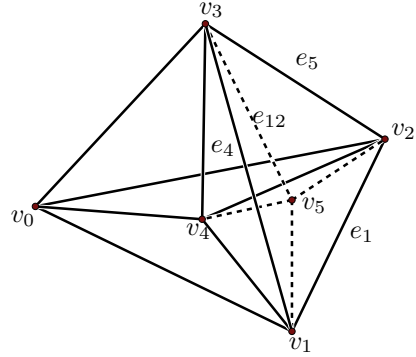}
    \label{fig:second-constr-2}}
  
  \subfigure[Part of the Schlegel diagram of $Q^{(3)}$
      where the edge $e_7$ (indicated as a dotted line)
      from $Q^{(2)}$ is no longer an edge.]{
    \psfrag{e0}[tr][tr]{$e_0$}
    \psfrag{e1}[tl][tl]{$e_1$}
    \psfrag{e2}[tr][tr]{$e_2$}
    \psfrag{e3}[br][br]{$e_3$}
    \psfrag{e4}[tr][tr]{$e_4$}
    \psfrag{e5}[bl][bl]{$e_5$}
    \psfrag{e6}[br][br]{$e_6$}
    \psfrag{e7}[tr][tr]{$e_7$}
    \psfrag{e8}[bl][bl]{$e_8$}
    \psfrag{e9}[br][br]{$e_9$}
    \psfrag{e16}[br][br]{$e_{16}$}
    \psfrag{v0}[t][t]{$v_0$}
    \psfrag{v1}[l][l]{$v_1$}
    \psfrag{v2}[bl][bl]{$v_2$}
    \psfrag{v3}[bl][bl]{$v_3$}
    \psfrag{v4}[bl][bl]{$v_4$}
    \psfrag{v5}[bl][bl]{$v_5$}
    \psfrag{v6}[bl][bl]{$v_6$}
    \includegraphics[width=.40\textwidth]{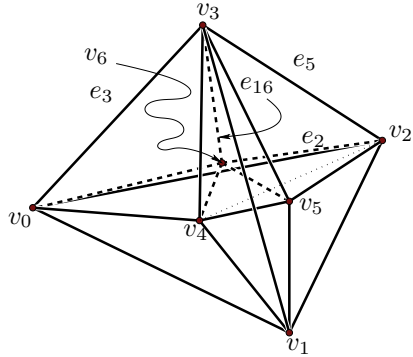}
    \label{fig:second-constr-3}}
\hspace{.05\textwidth}
  \subfigure[Part of the Schlegel diagram of $Q^{(4)}$
      with the simplex facet $F(v_4,v_5,v_6,v_7)$
      obtained in the last step.]{
    \psfrag{e10}[tr][tr]{$e_{10}$}
    \psfrag{e11}[tr][tr]{$e_{11}$}
    \psfrag{e12}[tr][tr]{$e_{12}$}
    \psfrag{e13}[tr][tr]{$e_{13}$}
    \psfrag{v0}[t][t]{$v_0$}
    \psfrag{v1}[l][l]{$v_1$}
    \psfrag{v2}[bl][bl]{$v_2$}
    \psfrag{v3}[bl][bl]{$v_3$}
    \psfrag{v4}[bl][bl]{$v_4$}
    \psfrag{v5}[bl][bl]{$v_5$}
    \psfrag{v6}[bl][bl]{$v_6$}
    \psfrag{v7}[bl][bl]{$v_7$}
    \includegraphics[width=.40\textwidth]{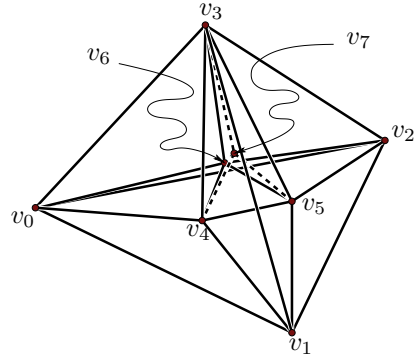}
    \label{fig:second-constr-4}}
\caption{The second construction.}
\end{figure}
\begin{table}[bt]
  \small
  \caption{Parameters describing the second construction.}
  \label{tab:second}
  \begin{sideways}
    \begin{tabular}{ccccccc}
      step $i$ & base facet $S_i$ & facets in ${\mathcal F}_i$ & facets
      in ${\mathcal N}_i$ & new facets & changing edges & degree \\
      \hline
      $1$                              & $S$                  & $F_S(e_0,e_3,e_4)$         & --                         & $F(v_0,v_1,v_2,v_4)$ & $[v_0,v_4]$         & $3$ \\
      (Fig.~\ref{fig:second-constr-1}) &                      &                            &                            & $F(v_1,v_2,v_3,v_4)$ & $[v_1,v_4]$         & $3$ \\
      \mbox{}                          &                      &                            &                            & $F(v_0,v_2,v_3,v_4)$ & $[v_2,v_4]$         & $3$ \\
      \mbox{}                          &                      &                            &                            &                      & $e_9:=[v_2,v_4]$    & $3$ \\
      \mbox{}                          &                      &                            &                            &                      & $e_1$               & $d_1+1$ \\
      \mbox{}                          &                      &                            &                            &                      & $e_2$               & $d_2+1$ \\
      \mbox{} & & & & & $e_5$ & $d_5+1$ \\ \hline
      $2$                              & $F(v_1,v_2,v_3,v_4)$ & $F(v_0,v_1,v_2,v_4)$       & --                         & $F(v_1,v_3,v_4,v_5)$ & $[v_1,v_5]$         & $3$ \\
      (Fig.~\ref{fig:second-constr-2}) &                      & $F_{S_2}(e_1,e_4,e_5)$     &                            & $F(v_2,v_3,v_4,v_5)$ & $[v_2,v_5]$         & $3$ \\
      \mbox{}                          &                      &                            &                            &                      & $e_{12}:=[v_3,v_5]$ & $3$ \\
      \mbox{}                          &                      &                            &                            &                      & $[v_4,v_5]$         & $3$ \\
      \mbox{}                          &                      &                            &                            &                      & $e_1$               & $d_1$ \\
      \mbox{} & & & & & $e_9$ & $4$ \\ \hline
      $3$                              & $F(v_0,v_2,v_3,v_4)$ & $F(v_0,v_1,v_2,v_4,v_5)$   & $F(v_2,v_3,v_4,v_5)$       & $F(v_0,v_3,v_4,v_6)$ & $[v_0,v_6]$         & $3$ \\
      (Fig.~\ref{fig:second-constr-3}) &                      & $F_{S_3}(e_2,e_3,e_5)$     &                            & $F(v_2,v_3,v_5,v_6)$ & $[v_2,v_6]$         & $3$ \\
      \mbox{}                          &                      &                            &                            & $F(v_3,v_4,v_5,v_6)$ & $e_{16}:=[v_3,v_6]$ & $4$ \\
      \mbox{}                          &                      &                            &                            &                      & $[v_4,v_6]$         & $3$ \\
      \mbox{}                          &                      &                            &                            &                      & $[v_5,v_6]$         & $3$ \\
      \mbox{}                          &                      &                            &                            &                      & $e_2$               & $d_2$ \\
      \mbox{}                          &                      &                            &                            &                      & $e_5$               & $d_5$ \\
      \mbox{} & & & & & $e_{12}$ & $4$ \\ \hline
      $4$                              & $F(v_3,v_4,v_5,v_6)$ & $F(v_0,v_3,v_4,v_6)$       & --                         & $F(v_4,v_5,v_6,v_7)$ & $[v_3,v_7]$         & $3$ \\
      (Fig.~\ref{fig:second-constr-4}) &                      & $F(v_1,v_3,v_4,v_5)$       &                            &                      & $[v_4,v_7]$         & $3$ \\
      \mbox{}                          &                      & $F(v_2,v_3,v_5,v_6)$       &                            &                      & $[v_5,v_7]$         & $3$ \\
      \mbox{}                          &                      &                            &                            &                      & $[v_6,v_7]$         & $3$ \\
      \mbox{}                          &                      &                            &                            &                      & $e_9$               & $3$ \\
      \mbox{}                          &                      &                            &                            &                      & $e_{12}$            & $3$ \\
      \mbox{} & & & & & $e_{16}$ & $3$ \\ \hline
    \end{tabular}
  \end{sideways}
\end{table}

\begin{definition}
  Let $Q$   be a $4$-polytope  with   a simplex facet $S$   in bounded
  position.  We denote the polytope $Q^{(4)}$ obtained by applying the
   four steps  of Table~\ref{tab:second} by $\second(Q;S)$.
\end{definition}
Remark~\ref{rem:ambiguous} also applies to $\second(P;S)$.
As in the first construction, we  have again created simplex facets in
bounded position in $\second(Q;S)$ and $2$-simplicity is preserved.
\begin{lemma}\label{lemma:bounded-2} 
  Let $Q$   be  a $4$-polytope with   a simplex  facet  $S$ in bounded
  position and the vertices of $S$ numbered  in arbitrary order.  Then
  the simplex facet of $\second(Q;S)$  constructed in the last step is
  in bounded position.  Furthermore, all  edges of $Q$ are again edges
  of  $\second(Q;S)$ with  the  same degree,  and  $\fdeg(e)=3$ for all
  edges $e$ in $\second(Q;S)$ that were not edges in $Q$.\qed
\end{lemma}

%%%%%%%%%%%%%%%%%%%%%%%%%%%%%%%%%%%%%%%%%%%%%%%%%%

\subsection{Properties of the Constructions}

The two constructions preserve the properties we are interested in.
\begin{theorem}\label{thm:f-vector-4} 
  Let  $P$ be  a  $4$-polytope with  a   simplex facet $S$  in bounded
  position. Then
    \[ f(\first(P;S)) \; = \; f(P) + (5,20,20,5) , \]
    and
    \[ f(\second(P;S)) \; = \; f(P) + (4,16,16,4) . \]
\end{theorem}
\begin{proof}
  In every step of both constructions, we add one vertex each by Corollary~\ref{cor:vertices-plus-1}.
  Hence $f_0(\first(P;S)) = f_0(P)+5$ and $f_0(\second(P;S)) = f_0(P)+4$.
  \par
  Propositions~\ref{prop:new-edges}  and   \ref{prop:new-edges-N} tell
  how many   edges are  added  in  each  step. Accordingly, the  first
  construction adds $4$ edges in each step. In Steps 1, 2 and 4 of the
  second construction, we also  add $4$ edges  each; in the third step
  $5$ edges  are  added, but  one edge  is destroyed,   so the overall
  change is also $4$. Therefore, we  have $f_1(\first(P;S)) = f_1(P) +
  5 \cdot 4 = f_1(P)+20$ and $f_1(\second(P;S)) = f_1(P) + 4 \cdot 4 =
  f_1(P)+16$. Finally, the  number of  new  facets in  the  respective
  steps   can    be   read  of    from    Tables~\ref{tab:first}   and
  \ref{tab:second}  (or Proposition~\ref{prop:facets-ps}) -- note that
  the facets $S_i$ in every step, as well as the facet of $Q^{(3)}$ in
  ${\mathcal N}_3$, disappear. Summarising, we get $f_3(\first(P;S)) =
  f_3(P)+5$ and $f_3(\second(P;S)) =  f_3(P)+4$.
  \par
  By Euler's equation, the number of ridges is determined by
  $f_0$, $f_1$ and $f_3$, which implies the claim.
\end{proof}

\begin{corollary} \label{cor:g2}
  Let $P$ be a $4$-polytope with a simplex facet $S$ in bounded position.
  Then $g_2(\first(P;S)) = g_2(\second(P;S)) = g_2(P)$.
\end{corollary}
\begin{proof}
  Define $a(P) := f_{02}(P) - 3f_2(P)$ and $b(P) := f_1(P) - 4f_0(P)$.
  Then $g_2(P) = a(P) + b(P) + 10$.
  \par
  We first show that $a(\first(P;S)) = a(P)$. Any $2$-face of $\first(P;S)$
  that is not a face of $P$ is a triangle by Proposition~\ref{prop:facets-ps}.
  Furthermore, any $2$-face of $P$ that is not a face of $\first(P;S)$
  is a ridge of the base facet in one of the five steps; since all base facets
  are simplices, all such $2$-faces are also triangles.
  Therefore, none of these faces causes a change in $a(P)$.
  \par
  For $b(P)$, we have, by Theorem~\ref{thm:f-vector-4},
  \begin{align*}
    b(\first(P;S)) & = f_1(\first(P;S)) - 4 f_0(\first(P;S)) \\
    & = f_1(P)+20 - 4(f_0(P)+5) \; = \; f_1(P) - 4 f_0(P) \; = b(P) .
  \end{align*}
  This implies $g_2(\first(P;S)) = g_2(P)$.
  \par
  Similar reasoning shows the claim for $\second(P;S)$.
\end{proof}

\begin{theorem}\label{thm:prop-2s2s}
  Let $P$ be a \tstsfp\ with a simplex facet  $S$ in bounded position.
  Then  $\first(P;S)$ and $\second(P;S)$   are  again $2$-simple   and
  $2$-simplicial.
\end{theorem}
\begin{proof}
  By   Theorem~\ref{thm:k-simplicial},   all intermediate polytopes in
  both constructions are $2$-simplicial,  hence also $\first(P;S)$ and
  $\second(P;S)$.   Furthermore,   all edges   in  $\first(P;S)$   and
  $\second(P;S)$  have degree $3$  by Lemmas~\ref{lemma:bounded-1} and
  \ref{lemma:bounded-2}, and since this is true for $P$.
\end{proof}

\begin{remark}
  There are more sequences of pseudo-stacking operations that preserve
  $g_2=0$ than the two given here.   Even more can be constructed with
  the   help of a further   generalised pseudo-stacking operation, see
  Remark~\ref{rem:generalised-pseudo}.

  A more symmetric version using these generalisations also applies in
  dimensions $d\ge 5$,  i.e.\ using this  construction  one can obtain
  elementary $2$-simple and $2$-simplicial   $d$-polytopes 
  with arbitrarily large numbers of vertices.
\end{remark}

%%%%%%%%%%%%%%%%%%%%%%%%%%%%%%%%%%%%%%%%%%%%%%%%%%%%%%%%%%%%%%%%%%%%%%%%%%%%%%%

\section{Examples and Results} \label{sec:ex}

Now we are ready to prove our main result: that $\ell_1$ is in fact an
extremal  ray of $\fhull$. We  prove this  by providing three examples
$P_9$, $P_{10}$, and $P_{11}$ of $2$-simple, $2$-simplicial elementary
$4$-polytopes  with $9$, $10$,  and  $11$  vertices that have  simplex
facets     in    bounded  position.   Using     these     as input for
Theorem~\ref{thm:f-vector-4} we obtain  such polytopes for arbitrarily
high numbers of vertices.

%%%%%%%%%%%%%%%%%%%%%%%%%%%%%%%%%%%%%%%%%%%%%%%%%%

\subsection{Examples}

For each of the following three examples we give their vertex-facet incidences 
together with a Schlegel diagram.
Explicit     rational     coordinates are   in Table~\ref{tab:coords}.
Calculations for the realisations were done by computer with the
{\tt polymake} system \cite{GawrilowJoswig2}. Clients producing the
examples below and applying the constructions,  as well as coordinates
for   many  more   elementary $2$-simple   and $2$-simplicial  $4$-polytopes 
are available from the authors.

The combinatorial descriptions of the examples, together with many 
more such polytopes having a small 
number of vertices, were found via a complete enumeration 
approach using a client for the {\tt polymake} system. 

\begin{figure}[bt]
  \begin{minipage}[b]{.4\textwidth}
    \begin{viftable}
      \{ v_1,v_2,v_3,v_4,v_5 \} \\
      \{ v_3,v_4,v_5,v_7,v_8 \} \\
      \{ v_5,v_6,v_7,v_8 \} \\
      \{ v_2,v_4,v_5,v_6,v_8 \} \\
      \{ v_0,v_1,v_2,v_3 \} \\
      \{ v_0,v_2,v_3,v_5,v_6,v_7 \} \\
      \{ v_0,v_4,v_6,v_7,v_8 \} \\
      \{ v_0,v_1,v_3,v_4,v_7 \} \\
      \{ v_0,v_1,v_2,v_4,v_6 \}
    \end{viftable}
  \end{minipage}
  \hfill
  \begin{minipage}[b]{.55\textwidth}
    \psfrag{v0}[tr][tr]{$v_0$}
    \psfrag{v1}[tr][tr]{$v_1$}
    \psfrag{v2}[tr][tr]{$v_2$}
    \psfrag{v3}[tr][tr]{$v_3$}
    \psfrag{v4}[tr][tr]{$v_4$}
    \psfrag{v5}[tr][tr]{$v_5$}
    \psfrag{v6}[tr][tr]{$v_6$}
    \psfrag{v7}[tr][tr]{$v_7$}
    \psfrag{v8}[tr][tr]{$v_8$}
    \includegraphics[width=.85\textwidth]{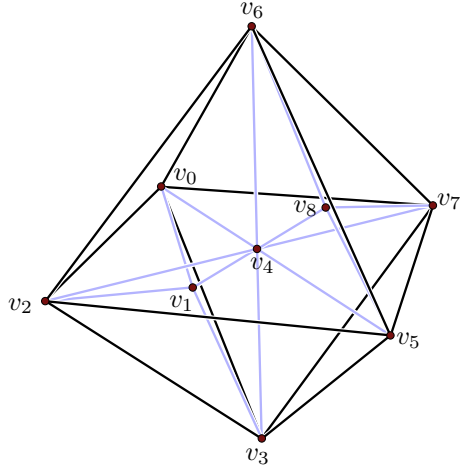}
  \end{minipage}
  \centering
  \caption{The vertex-facet-incidences and a Schlegel diagram of the polytope $P_9$.}
  \label{fig:w9}
\end{figure}
\begin{compactenum}
\item The smallest non-trivial \tstsfp\ is $P_9$ (see also the remark below).
  See Figure~\ref{fig:w9}.
  \begin{figure}[bt]
    \begin{minipage}[t]{.36\textwidth}
      \begin{viftable}
        \{ v_3,v_4,v_7,v_8,v_9 \} \\
        \{ v_1,v_2,v_3,v_4,v_6,v_9 \} \\
        \{ v_1,v_2,v_4,v_5 \} \\
        \{ v_1,v_3,v_7,v_9 \} \\
        \{ v_0,v_2,v_3,v_6 \} \\
        \{ v_0,v_4,v_7,v_8 \} \\
        \{ v_0,v_1,v_4,v_5,v_7,v_9 \} \\
        \{ v_0,v_1,v_2,v_5,v_6 \} \\
        \{ v_0,v_2,v_3,v_4,v_5,v_8 \} \\
        \{ v_0,v_1,v_3,v_6,v_7,v_8 \}
      \end{viftable}
    \end{minipage}
    \hfill
    \begin{minipage}[t]{.62\textwidth}
      \vspace{.5cm}
      \hspace{.5cm}
      \psfrag{v0}[tr][tr]{$v_0$}
      \psfrag{v1}[tr][tr]{$v_1$} \psfrag{v2}[tr][tr]{$v_2$}
      \psfrag{v3}[tr][tr]{$v_3$} \psfrag{v4}[tr][tr]{$v_4$}
      \psfrag{v5}[tr][tr]{$v_5$} \psfrag{v6}[tr][tr]{$v_6$}
      \psfrag{v7}[tr][tr]{$v_7$} \psfrag{v8}[tr][tr]{$v_8$}
      \psfrag{v9}[tr][tr]{$v_9$}
      \hspace{-.1\textwidth}\includegraphics[width=.99\textwidth]{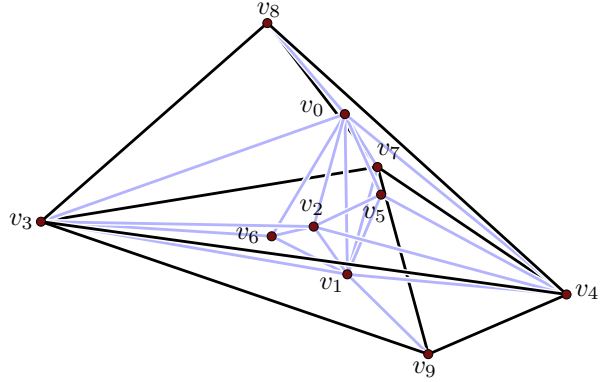}
    \end{minipage}
    \centering
    \caption{The vertex-facet-incidences and a Schlegel diagram of the polytope $P_{10}$.}
    \label{fig:p10}
  \end{figure}
\item The second example, $P_{10}$,  has $10$ vertices.
  See Figure~\ref{fig:p10}.
  Note that   $P_{10}$   is  not  combinatorially  equivalent   to the
  $4$-dimensional hypersimplex, nor its  dual --  although it has  the
  same  flag vector.  Nevertheless,  the   hypersimplex itself can  be
  obtained    with    the    extended  construction    described    in
  Remark~\ref{rem:generalised-pseudo}.
  
  \begin{figure}[tb]
    \begin{minipage}[t]{.4\textwidth}
      \begin{viftable}
        \{ v_1,v_5,v_6,v_7,v_9 \} \\
        \{ v_2,v_3,v_4,v_7,v_8,v_{10} \} \\
        \{ v_3,v_4,v_5,v_{10} \} \\
        \{ v_3,v_5,v_6,v_7,v_{10} \} \\
        \{ v_1,v_2,v_3,v_8 \} \\
        \{ v_1,v_3,v_6,v_7,v_8 \} \\
        \{ v_0,v_2,v_4,v_7 \} \\
        \{ v_0,v_1,v_5,v_9 \} \\
        \{ v_0,v_1,v_2,v_3,v_4,v_5,v_6 \} \\
        \{ v_0,v_4,v_5,v_7,v_9,v_{10} \} \\
        \{ v_0,v_1,v_2,v_7,v_8,v_9 \}
      \end{viftable}
    \end{minipage}
    \hfill
    \begin{minipage}[t]{.55\textwidth}
      \vspace{.3cm}
      \psfrag{v0}[tr][tr]{$v_0$}
      \psfrag{v1}[tr][tr]{$v_1$} \psfrag{v2}[tr][tr]{$v_2$}
      \psfrag{v3}[tr][tr]{$v_3$} \psfrag{v4}[tr][tr]{$v_4$}
      \psfrag{v5}[tr][tr]{$v_5$} \psfrag{v6}[tr][tr]{$v_6$}
      \psfrag{v7}[tr][tr]{$v_7$} \psfrag{v8}[tr][tr]{$v_8$}
      \psfrag{v9}[tr][tr]{$v_9$} \psfrag{v10}[tr][tr]{$v_{10}$}
      \hspace{-.1\textwidth}\includegraphics[width=.99\textwidth]{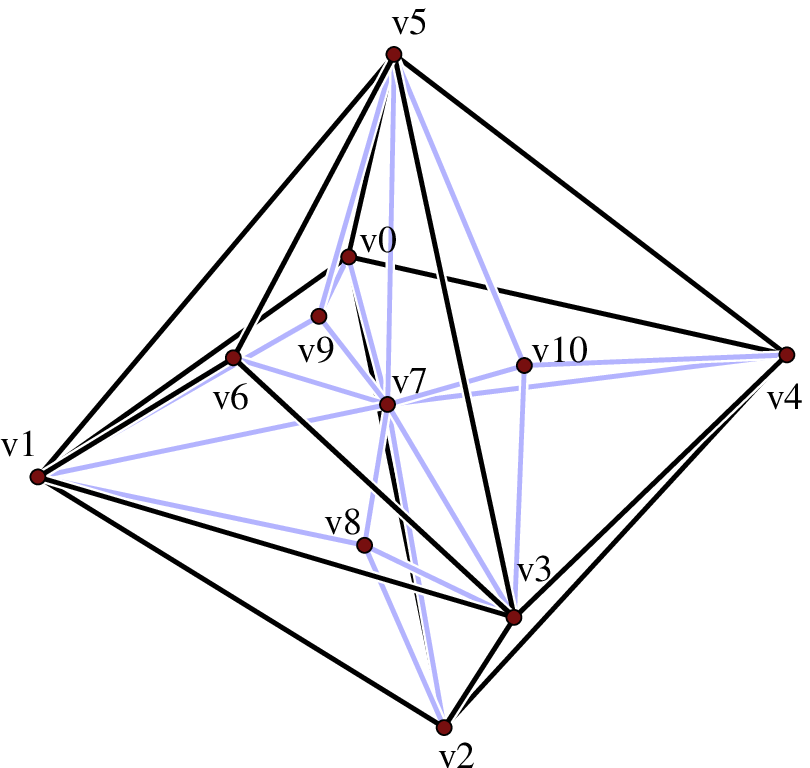}
    \end{minipage}
    \centering
    \caption{The vertex-facet-incidences and a Schlegel diagram of the polytope $P_{11}$.}
    \label{fig:p11}
  \end{figure}
\item The last example, $P_{11}$, has $11$  vertices. 
  See      Figure~\ref{fig:p11}. 
\end{compactenum}
\begin{remark}
  $P_9$ and $P_{10}$ are both self-dual. $P_9$ is the unique non-trivial
  $2$-simple and $2$-simplicial $4$-polytope with smallest number of vertices,
  i.e.\ except for the simplex $\Delta_4$ there are no further
  \tstsfp{s} with less than $10$ vertices. We do not prove this here; 
  part of the proof can be found in \cite{DissA}.
\end{remark}

\begin{remark}
  The  two  examples $P_9$ and  $P_{10}$ can  in fact  themselves be
  obtained by  applying the  two constructions to   the $4$-simplex:
  $P_9 = \second(\Delta_4;S)$ and $P_{10} = \first(\Delta_4;S)$ with
  an   arbitrary facet $S$ of  $\Delta_4$.    Note however, that the
  simplex  has  no  facet  in bounded  position;  nevertheless,  the
  constructions  remain valid, since   in both cases the  first step
  requires the  added vertex  to  lie in  only  one adjacent   facet
  hyperplane.

  Also, $P_{11}$ can be  obtained  via three pseudo-stacking  steps in
  the following way.  Take a $3$-dimensional octahedron $O^3$.  Choose
  a facet $R$ of $O^3$  and let $R_1$, $R_2$, and  $R_3$ be the facets
  adjacent to  $R$.  Let $B^3$ be the  polytope obtained from $O^3$ by
  stacking the  facet $R$.  Let $PB^3$  be the pyramid over $B^3$ with
  apex $v$ and $F_1$, $F_2$, and $F_3$ be the facets of $PB^3$ arising
  as  pyramids over $R_1$, $R_2$, and  $R_3$.  $PB^3$ is an elementary
  $2$-simplicial polytope, but it is not $2$-simple.

  Let  $\mathcal   F_i$ be the   sets  of   facets  adjacent to  $F_i$
  containing   the   vertex  $v$.  Then    $|\mathcal   F_i|=3$.   Set
  $P^{(0)}:=PB^3$ and define polytopes $P^{(i)}:=\pstack{F_i}{\mathcal
    F_i}{P^{(i-1)}}$ for $i=1,2,3$.  Then $P_{11}=P^{(3)}$.  Note that
  some facets in $\mathcal F_2$ and $\mathcal  F_3$ are stacked by the
  previous steps, but this does not influence the result.

  $P_{11}$  demonstrates that one  can obtain elementary $2$-simple and
  $2$-simplicial   $4$-polytopes by  pseudo-stacking polytopes without
  these properties in a suitable way.
\end{remark}
\begin{table}[bt]
  \small
  \centering
  \caption{Coordinates for the three examples  of
    $2$-simple and $2$-simplicial  $4$-polytopes that we have
    discussed in this section.}
  \begin{tabular}{lll}
    \begin{minipage}[t]{.27\textwidth}
      \centering
      $P_9$\\      
      \begin{verticestable}
        3 &  0 &  0 & 0 \\
        1 &  1 &  1 & 1 \\
        0 &  3 &  0 & 0 \\
        0 &  0 &  3 & 0 \\
        0 &  0 &  0 & \frac{3}{2} \\
        -3 &  0 &  0 & 0 \\
        0 &  0 & -3 & 0 \\
        0 & -3 &  0 & 0 \\
        -1 & -1 & -1 & 1
      \end{verticestable}
    \end{minipage}
    &
    \begin{minipage}[t]{.32\textwidth}
      \centering
      $P_{10}$\\      
      \begin{verticestable}
        9 & -3 & -3 & -3 \\
        -3 &  9 & -3 & -3 \\
        -3 & -3 & -3 & -3 \\
        -3 & -3 &  9 & -3 \\
        -3 & -3 & -3 &  9 \\
        1 & -3 & -7 &  1 \\
        -3 &  1 &  1 & -7 \\
        3 &  3 &  3 &  3 \\
        5 & -3 &  5 &  1 \\
        -3 &  5 &  1 &  5
      \end{verticestable}
    \end{minipage}
    &
    \begin{minipage}[t]{.32\textwidth}
      \centering
      $P_{11}$\\      
      \begin{verticestable}
        1 & 0 & 0 & 0 \\
        0 & 1 & 0 & 0 \\
        0 & 0 & 1 & 0 \\
        -1 & 0 & 0 & 0 \\
        0 & -1 & 0 & 0 \\
        0 & 0 & -1 & 0 \\
        -\frac{11}{21} & \frac{11}{21} & -\frac{11}{21} & 0 \\
        -\frac{11}{147} & \frac{11}{147} & -\frac{11}{147} & 1 \\
        -\frac{428}{1617} & \frac{428}{1617} & \frac{68}{147} & \frac{1}{2} \\
        \frac{68}{147} & \frac{428}{1617} & -\frac{428}{1617} & \frac{1}{2} \\
        -\frac{428}{1617} & -\frac{68}{147} & -\frac{428}{1617} & \frac{1}{2}
      \end{verticestable}
    \end{minipage}
  \end{tabular}
  \label{tab:coords}
\end{table}

%%%%%%%%%%%%%%%%%%%%%%%%%%%%%%%%%%%%%%%%%%%%%%%%%%

\subsection{Conclusions}

With     these  examples        and     the     constructions       of
Section~\ref{sec:constr-2s2s}, $2$-simplicial and $2$-simple $4$-polytopes
with $g_2=0$ can be found for almost all numbers of vertices. 
Here is the main theorem of this paper.
\begin{theorem}\label{thm:many-polys}
  Elementary \tstsfp{s} with $k$  vertices exist for $k=5,9,10,11$ and
  $k \geq 13$.
\end{theorem}
\begin{proof}
  We show that for $k$ as given in the claim there exist elementary \tstsfp{s}
  with $k$ vertices that have at least one simplex facet in bounded position.
  \par
  For $k=9,10,11$, the above examples $P_9$, $P_{10}$, resp.\ $P_{11}$
  have the desired properties.
  Let $k \geq 13$. By induction there is such a polytope $P$ with
  $k-5$ or $k-4$ vertices and a simplex facet $S$ in bounded position.
  Then the polytope $\first(P;S)$, resp.\ $\second(P;S)$ are elementary \tstsfp{s}
  on $k$ vertices by Lemma~\ref{lemma:bounded-1} resp.~\ref{lemma:bounded-2}
  and Corollary~\ref{cor:g2}.
\end{proof}
\begin{corollary}\label{cor:ray} 
  The ray $\ell_1$ is contained in the convex hull of all flag vectors
  of $4$-polytopes.\qed
\end{corollary}
There are also many \tstsfp{s} with $g_2>0$; still the
existence of a \tstsfp\ with $12$ vertices is an open question.

%%%%%%%%%%%%%%%%%%%%%%%%%%%%%%%%%%%%%%%%%%%%%%%%%%%%%%%%%%%%%%%%%%%%%%%%%%%%%%%

\providecommand{\bysame}{\leavevmode\hbox to3em{\hrulefill}\thinspace}
\providecommand{\MR}{\relax\ifhmode\unskip\space\fi MR }
\providecommand{\MRhref}[2]{%
  \href{http://www.ams.org/mathscinet-getitem?mr=#1}{#2}
}
\providecommand{\href}[2]{#2}

\end{document}